\documentclass[a4paper,12pt]{article}

\usepackage{amssymb,amsthm,amscd}
\usepackage{amsmath}
\usepackage{subfigure}
\usepackage{graphicx}
\usepackage{color}

\usepackage{enumitem}
\setitemize{noitemsep,leftmargin=15pt, topsep=0pt,parsep=0pt,partopsep=10pt}
\setenumerate{noitemsep,leftmargin=25pt, topsep=0pt,parsep=0pt,partopsep=0pt}

\theoremstyle{plain}
\newtheorem{theorem}{Theorem}[section]

\newtheorem{prop}[theorem]{Proposition}
\newtheorem{la}[theorem]{Lemma}
\newtheorem{cor}[theorem]{Corollary}
\theoremstyle{remark}
     \newtheorem{remark}[theorem]{Remark}

\def\barr{\begin{array}}
\def\earr{\end{array}}
\def\beqarr*{\begin{eqnarray*}}
\def\eeqarr*{\end{eqnarray*}}

\def\mapright#1{\smash{\mathop{\longrightarrow}\limits^{#1}}}
\def\maprdown#1{\Big\downarrow
                   \rlap{$\vcenter{\hbox{$\scriptstyle#1$}}$}}
\def\mapldown#1{\llap{$\vcenter{\hbox{$\scriptstyle#1$}}$}\Big\downarrow}

\newcommand{\N}{\mathbb{N}}

\newcommand{\R}{\mathbb{R}}

\newcommand{\Z}{\mathbb{Z}}

\begin{document}

\pagenumbering{arabic}


\title{\Large\bf A Proof of the Persistence of Anti-integrable States for Three-Dimensional Quadratic Diffeomorphisms}

\author{Yi-Chiuan Chen \\
  Institute of  Mathematics, Academia Sinica, Taipei 106319, Taiwan\\
   Email: YCChen@math.sinica.edu.tw \\}

\date{December 5, 2024}

\maketitle

\begin{abstract}
  Three-dimensional quadratic diffeomorphisms with quadratic inverse generically have five independent parameters. When some parameters approach infinity, the diffeomorphisms may exhibit a so-called anti-integrable limit
 in the traditional sense of Aubry and Abramovici. That is, the dynamics of the diffeomorphisms reduce to symbolic dynamics on finite number of symbols. However, the diffeomorphisms may reduce to quadratic correspondences when parameters approach infinity, and the traditional anti-integrable limit does not deal with this situation. Meiss asked what about an anti-integrable limit for it. A remarkable progress was achieved very recently by the work of Hampton and Meiss [{\it SIAM J. Appl. Dyn. Syst.} 21 (2022), pp. 650--675]. Using the contraction mapping theorem, they showed there is a bijection between the anti-integrable states and the sequences of branches of a quadratic correspondence. They also showed that an anti-integrable state  can be continued to a genuine orbit of the three-dimensional diffeomorphism.  This paper aims to contribute the progress, by means of the implicit function theorem. We shall show that, under a slightly more restricted condition than that imposed by Hampton and Meiss, the bijection indeed is a topological conjugacy and establish the uniform hyperbolicity of the continued genuine orbits. 
\end{abstract}





\tableofcontents

\section{Introduction}

It is shown in \cite{LLM1999, LM1998} that any  three-dimensional  quadratic diffeomorphism with quadratic inverse can be written (up to an affine transformation) in the following form 
\begin{equation}
(x,y,z)\mapsto (\alpha+\tau x-\sigma y +\delta z+ax^2+bxy+cy^2, x, y),  \label{map0}
\end{equation}
characterized by seven parameters,  in which the parameter 
$\delta$ is the Jacobian determinant.  
Providing $a+b+c\not= 0$ and $2a+b\not=0$, there exists an affine coordinate transformation with which the above diffeomorphism takes the  form 
\begin{equation}
 L: (x,y,z)\mapsto (\alpha-\sigma y +\delta z+ax^2+bxy+cy^2, x, y) \label{map1}
\end{equation}
(namely $\tau=0$) subject to 
\[
   a+b+c=1. 
\]
 If $a+b+c=0$  or $2a+b=0$, then there are other transformations with which two of the seven parameters in \eqref{map0} can be eliminated.
 The map $L$ given by the diffeomorphism \eqref{map1}  has five independent parameters and is what we shall consider in this paper. 

We shall choose $\alpha$, $\sigma$, $\delta$, $a$  and $c$ to be the five parameters, and are concerned with the dynamics of $L$ when $\alpha$, $\sigma$ or $\delta$ becomes infinite. 
(The same as Hampton and Meiss \cite{HM2022, HM2024} did, we assume $a$ and $c$ are ``structural" parameters that remain finite.) 

Let $(x_t,y_{t}, z_t)_{t\in\mathbb{Z}}$ with $(x_t,y_t, z_t)=L^t(x,y,z)$. 
Then, 
 $(x_t,y_{t}, z_t)_{t\in\mathbb{Z}}$ is an orbit of $L$ if and only if the sequence $(x_t)_{t\in\mathbb{Z}}$ is a solution of the following third-order difference equation
\[
  x_{t+1}=\alpha-\sigma x_{t-1}+\delta x_{t-2}+a x_t^2+bx_t x_{t-1}+c x_{t-1}^2, \quad \forall t\in\mathbb{Z},  
\]
  upon noting that $y_{t+1}=x_t$, $z_{t+1}=x_{t-1}$.
If we define
 $\xi_t=\epsilon x_t$ for $\epsilon>0$, then the difference equation above becomes 
\begin{equation}
 \epsilon^2 \alpha -\epsilon (\xi_{t+1}+\sigma \xi_{t-1}-\delta \xi_{t-2})+a\xi_t^2+b\xi_{t}\xi_{t-1}+c\xi_{t-1}^2=0. \label{Recurrence}
\end{equation} 
In this paper, the {\it anti-integrable (AI) limit} is the limit $\epsilon\to 0$. The limit is singular in the sense that the third-order difference equation \eqref{Recurrence} degenerates to a lower-order one or to an algebraic equation.  

How fast  the parameters $\alpha$, $\sigma$ and $\delta$ approach to infinite will determine the dynamics of \eqref{Recurrence} at the AI limit.
Suppose that the approach of parameters to infinity is characterized by $\epsilon$: Assume that
$\alpha$, $\sigma$ and $\delta$ depend $C^1$ on $\epsilon$,  and that
\[  \lim_{\epsilon\to 0} \epsilon^2 \alpha = \alpha_1, \qquad
 \lim_{\epsilon\to 0} \epsilon \sigma =\sigma_1, \qquad
 \lim_{\epsilon\to 0} \epsilon \delta =\delta_1. 
\]
 For instance, $\alpha_1=0$ if $\alpha=1/\epsilon$, whereas $\alpha_1=1$ if $\alpha=1/\epsilon^2$. 
Note that if $\alpha_1$, $\sigma_1$ or $\delta_1$ is finite and not equal to zero, then when $\epsilon$ is small, $\alpha$, $\sigma$ or $\delta$ is of $O(\epsilon^{-2})$, $O(\epsilon^{-1})$ and  $O(\epsilon^{-1})$, respectively.  
Additional to the dependence on $\epsilon$, we also assume that $\epsilon^2\alpha$, $\epsilon \sigma$ and $\epsilon\delta$ depend jointly $C^1$  on $\alpha_1$,
 $\sigma_1$ and $\delta_1$, respectively. For instance, $\epsilon^2 \alpha\equiv\alpha_1$,
 $\epsilon\sigma=\sigma_1+\epsilon$, $\epsilon\delta=\delta_1+\epsilon$ or $\epsilon\delta=\delta_1-\epsilon^{3/2}\log\epsilon$. Hence\footnote{
For parameters close to infinity, it is convenient to introduce a factor $\epsilon$ and  rescale the phase and parameter spaces in term of $\epsilon$.
}, near and at the AI limit, the dynamics is determined implicitly by the following set of parameters $\epsilon$, $\alpha_1$, $\sigma_1$, $\delta_1$, $a$ and $c$.

In a special situation when $\delta=0$ (thus the map $L$ is essentially two-dimensional),  $a=1$, $\alpha=-1/\epsilon^2$ and  $\sigma$ is a constant (thus $\alpha_1=-1$, $\sigma_1=0$, $\delta_1=0$, $c=0$), the difference equation  \eqref{Recurrence} degenerates to the algebraic equation $\xi^2_t=1$ at the AI limit.
 The solution of the equation at the AI limit is an arbitrary sequence $(\xi_t)_{t\in\Z}\in\{-1,1\}^\Z$. 
Such
 sequences can be continued to solutions of \eqref{Recurrence} for $\epsilon$ sufficiently close to zero (see e.g. \cite{Aubr1995, SDM1999, SM1998}), hence give rise to genuine orbits of the H\'{e}non map $(x, y)\mapsto (-1/\epsilon^2-\sigma y+x^2, x)$ (or more precisely, of an ``embedded" H\'{e}non map $(x, y, z )\mapsto (-1/\epsilon^2-\sigma y+x^2, x, y)$ in the current setting).

The special situation above nicely demonstrate the concept of AI limit. This concept  was first introduced in 1990 by Aubry and  Abramovici \cite{AA1990} for  the standard map, and    has been applied to various classical dynamical systems, for instance,  one-dimensional polynomial maps \cite{Chen2007, Chen2008}, the H\'{e}non maps \cite{Chen2018, SDM1999, SM1998} (see also \cite{Qin2001} for high-dimensional H\'{e}non-like maps), the Smale horseshoe \cite{Chen2006} and high-dimensional symplectic maps \cite{Chen2005, MM1992}. 
 Briefly, the ``dynamics" at the limit  becomes non-deterministic, reducing to a Bernoulli shift or subshift of finite type on symbols. Any sequence allowed by the shift on the symbols is called an {\it AI orbit} or {\it AI state} (like the sequence $(\xi_t)_{t \in\Z} \in\{-1,1\}^\Z$ in the above paragraph). Remarkably, with some nondegeneracy conditions at the AI limit, every AI state persists to a genuine orbit of the original map. In addition, the collection of  initial points of these genuine orbits forms a chaotic invariant set.
The analysis of the persistence of AI orbits can be based on the contraction mapping theorem (CMT) \cite{AA1990, SM1998}, the implicit function  theorem (IFT) \cite{Aubr1995,Chen2006,MM1992}, or  numerical computation \cite{HM2022, SDM1999}. 

This concept and methodology of anti-integrability have shown powerful for proving existence of chaotic orbits not only in discrete-time systems (i.e. maps or difference equations) but also in continuous-time ones (e.g. Hamiltonian and Lagrangian systems \cite{BCM2013, BM1997, Chen2003}, $N$-body problems \cite{BM2000}, and billiards \cite{Chen2004}).  The readers are referred to \cite{Aubr1995, BT2015} for nice reviews. 

Now, let us return back to the difference equation \eqref{Recurrence}, Meiss (see for example \cite{Meis2015}) asked what about an AI limit for it. In general, the limit $\epsilon\to 0$ is quite different from the traditional AI limit: when $\epsilon\to 0$ (together with $\alpha\to\pm\infty$ and/or $\sigma\to\pm\infty$ but with $\delta$ remaining finite) the difference equation reduces to a quadratic correspondence, which is much more complicated than the algebraic equation, e.g. $\xi_t^2=1$, in the traditional limiting situation. 
Although some aspects of this question had been studied in \cite{JLM2008, LM2010}, where \eqref{Recurrence}  degenerates to a one-dimensional map   (see Subsection \ref{subsec:JLM} of  this paper),  the general case of the question was left very much open. (Despite that the one-dimensional map studied in \cite{JLM2008} is a fixed branch of a quadratic correspondence, it is deterministic, but the correspondence itself is not.)

A breakthrough was achieved very recently by Hampton and Meiss. In \cite{HM2022, HM2024}, they found parameter regions such that the quadratic correspondence at the AI limit introduces symbolic dynamics on two symbols (see Subsection \ref{subsec:HM}). The two  symbols represent two different branches, a pair of one-dimensional maps, of the correspondence. Each symbol sequence determines a sequence of one-dimensional maps whose orbits correspond to the AI states.    Moreover, each such AI state is unique and can be continued away from the AI limit, by means of a contraction operator,  becoming a genuine orbit of the three-dimensional map $L$. Their results generalize  but still capture the spirit of the original work of Aubry and Abramovici \cite{AA1990} on AI limit: the ``dynamics" at the AI limit is not deterministic. The paper \cite{HM2022} was a beautiful step towards the general AI limit of three (and possibly higher) dimensional maps.

The main purpose of this paper is to present a unified treatment of the aforementioned works \cite{HM2022, HM2024, JLM2008, LM2010} by making use of the  IFT. Although finding AI states is not a trivial task, our main concern is focusing on the persistence of AI states. Our approach relies on a notion of (uniform) hyperbolicity for a backward and forward  complete $C^1$ relation on $\mathbb{R}$ (see definition in Section \ref{sec:correspondence}). This notion is inspired by the work of \cite{HM2022}.

A byproduct of our approach is to confirm that those genuine orbits continued form the AI states in \cite{HM2022, HM2024} constitute a horseshoe. Hampton and Meiss's numerical results indicate that the set consisting of those orbit points has a Cantor-like structure. We shall see in Corollary \ref{cor:HM_horseshoe} that this set is indeed an invariant Cantor set and is hyperbolic for the map $L$ if a slightly more restricted condition than that in \cite{HM2022, HM2024} is imposed (see Remark \ref{condition_imposed}).  

Our main theorem (Theorem \ref{mainthm}), which states that all hyperbolic AI states persist to genuine orbits of the map $L$, is presented in Section \ref{sec:mainthm}. 
In Subsections  \ref{subsec:JLM} and \ref{subsec:HM}, we discuss and examine known examples of AI states in the literature and show that every example gives rise to a compact, hyperbolic, backward and forward complete set for the relation at the AI limit. 
In Section \ref{sec:discussion}, we discuss relevant issues for the AI limit. We first discuss an AI limit where the map $L$ reduces to a non-invertible 2D quadratic map rather than a correspondence (see Subsection \ref{subsec:Mira}). In Subsection \ref{subsec:equiva}, we provide an alternative proof of Lemma \ref{mainla}. A notion called {\it nondegeneracy} is defined in Subsection \ref{equiva2}. We shall see that this is an equivalent notion to  the hyperbolicity for a $C^1$ relation. 
We demonstrate in Subsection \ref{subsec:degenerate} a non-hyperbolic  example  in which  only two AI states can  be continued to become genuine orbits of the map $L$, but all others cannot. 
In the final section, we conclude with a summary of what are done in this paper.

\section{Quadratic correspondences and AI states} \label{sec:correspondence}

Let\footnote{Although we allow that $\delta$ grows to infinity as $\epsilon$ tends to zero, we shall treat in this paper the case $\delta=o(\epsilon^{-1})$ only so that $\epsilon\delta\to 0$ as $\epsilon\to 0$. Hence, for simplicity, the reader may think of $\delta$ as $O(1)$ or as a constant (e.g. \cite{HM2022, HM2024}). Nonetheless, the notion $\delta_1$ has the flexibility to deal with the case $\delta_1\not=0$ in future research. See Subsection \ref{subsec:Mira} for a discussion when ${\delta}_1$ is equal to some non-zero constant.}
 $e:=(\epsilon, \alpha_1, \sigma_1, \delta_1, a, c)=(\epsilon, \alpha_1, \sigma_1, 0, a, c)$ and  
\[ e^\dag:=(0, \bar{\alpha}_1, \bar{\sigma}_1, 0, \bar{a}, \bar{c}).
\] 
Let $\Xi$ be the Banach space $l_\infty(\mathbb{Z}, \mathbb{R}):=\{  \boldsymbol\xi |~\boldsymbol\xi=(\ldots,\xi_{-1},\xi_0,\xi_1,\ldots),~\xi_t\in\mathbb{R}~\forall t\in\Z,\ \mbox{bounded}\}$ of bounded sequences with the supremum 
norm $\|\cdot\|_\infty$.
Define a map $\mathcal{L}(\cdot; e):\Xi \to \Xi$,  $\mathcal{L}(\boldsymbol\xi; e)=(\mathcal{L}(\boldsymbol\xi; e)_t)_{t\in\Z}$,  by
\[ \mathcal{L}(\boldsymbol\xi; e)_t=\epsilon^2\alpha -\epsilon \xi_{t+1}-\epsilon\sigma \xi_{t-1}+\epsilon\delta \xi_{t-2}+a\xi_t^2+b\xi_{t}\xi_{t-1}+c\xi_{t-1}^2.
\]
Now, 
\begin{equation}
 \mathcal{L}(\boldsymbol\xi; e^\dag)_t=\bar{\alpha}_1 -\bar{\sigma}_1 \xi_{t-1}+{\bar a}\xi_t^2+\bar{b}\xi_{t}\xi_{t-1}+\bar{c}\xi_{t-1}^2, \label{L0delta2=0}
\end{equation}
with $\bar{b}=1-\bar{a}-\bar{c}$.

Observe that  a point $\boldsymbol{\xi}$ in $\Xi$ for which $\mathcal{L}(\boldsymbol\xi; e^\dag)=0$ is equivalent to a bounded ``orbit" $(\xi_t)_{t\in\Z}$ of the  following quadratic correspondence \cite{Bullet1988, McGe1992}
\begin{equation}
   \bar{\alpha}_1 -\bar{\sigma}_1 u+\bar{a}v^2+\bar{b}vu+\bar{c}u^2=0. \label{Q}
\end{equation}  
Notice that the ``dynamics" is no longer always deterministic. This correspondence can be viewed as a quadratic relation on $\R$
\begin{eqnarray}
  \mathcal{E} &:=& \mathcal{E}(\bar{\alpha}_1, \bar{\sigma}_1, \bar{a}, \bar{c}) \nonumber \\
  &=& \{(u,v):  \bar{\alpha}_1 -\bar{\sigma}_1 u+\bar{a}v^2+\bar{b}vu+\bar{c}u^2=0\}\subset \R^2. \label{relation}
\end{eqnarray}  
Recall that a {\it relation} on $\R$ is a subset of $\R\times \R$. A quadratic relation $\mathcal{E}$ is a conic (or called quadratic curve) on $\R\times\R$. If $\mathcal{E}$ is not a degenerate conic, depending on the value of the discriminant $\Delta\equiv \bar{b}^2-4\bar{a}\bar{c}$ being negative, equal to zero, or positive, $\mathcal{E}$ is either an ellipse,  a parabola,  or a hyperbola.

 Viewing a relation in terms of iteration, an iterate of $u$ under relation $\mathcal{E}$ is a point $v$ 
such that $(u, v)\in\mathcal{E}$. Such a point $v$ is not necessarily unique, nor does it necessarily exist.
A (bi-infinite) {\it orbit} of a relation $\mathcal{E}$ is a sequence $(\xi_t)_{t\in\Z}$ such that $(\xi_t, \xi_{t+1})\in\mathcal{E}$ for every $t$.

A subset of a relation is called $C^1$ if it is a $C^1$ embedded manifold of $\R\times \R$. The relation $\mathcal{E}$ defined by $\eqref{relation}$ certainly is $C^1$ for nondegenerate case. The {\it tangent relation} $T\mathcal{E}$ on $\R^2$ is the tangent bundle of a $C^1$ relation $\mathcal{E}$. For $(u, v)\in \mathcal{E}$, let $m=m(u,v)$ be the slope of $T_{(u,v)}\mathcal{E}$ viewed as a graph of a linear function of $\R$.  Direct calculation shows that 
\[
  m(u,v)=\frac{\bar{\sigma}_1-\bar{b}v-2\bar{c}u}{2\bar{a}v+\bar{b}u}.
\]

A compact set $\Lambda\subset \R$ is called a {\it (uniformly) expanding}  set for a $C^1$ relation $\mathcal{E}$ if there are constants  $C>0$ and $\lambda>1$ such that for every $\xi_0\in\Lambda$ the slope $m_t:=m(\xi_t, \xi_{t+1})$ of the linear relation $T_{(\xi_t, \xi_{t+1})}\mathcal{E}$ for the orbit $(\xi_t)_{t\in\Z}$ (whenever exists) with $(\xi_t, \xi_{t+1})\in\mathcal{E}\cap(\Lambda\times\Lambda)$  is not infinite and satisfies
\begin{equation}
 \prod_{n= 0}^{k-1} |m_{t+n}|\ge C\lambda^k~\quad \forall t\in \Z~\mbox{and}~ \forall k\ge 1. \label{Clambda+}
\end{equation}
It is called {\it (uniformly) contracting} if the slope $m_t$ is not zero and 
\begin{equation}
 \prod_{n= 0}^{k-1} |m_{t+n}|\le C^{-1}\lambda^{-k}~\quad \forall t\in \Z~ \mbox{and}~ \forall k\ge 1 . \label{Clambda-}
\end{equation}
We call a compact set $\Lambda\subset \R$ {\it (uniformly) hyperbolic} for a  relation if it is either expanding or contracting.  
We refer the reader to \cite{Sand1999} for more details about hyperbolicity for relations.
 
It is known that  an expanding compact invariant set $\Lambda$ of a $C^1$ interval map $f$ is {\it positively expansive}, that is,  there exists $\chi>0$ such that for distinct $x$, $\tilde{x}\in\Lambda$, there is some $n\ge 0$ such that $|f^n(x)-f^n(\tilde{x})|\ge\chi$.  We call a subset $\Lambda$ of $\R$ {\it expansive} for the quadratic correspondence \eqref{Q} or relation $\mathcal{E}$ if there exists $\chi>0$ such that for distinct $\xi_0$, $\tilde{\xi}_0\in\Lambda$ any orbits $(\xi_t)_{t\in\Z}$,  $(\tilde{\xi}_t)_{t\in\Z}$ satisfy $|\xi_n-\tilde{\xi}_n|\ge\chi$ for some $n\in\Z$. We have
\begin{prop} \label{Expansive}
A hyperbolic compact set for the  relation $\mathcal{E}$ is expansive.
\end{prop} 
\proof
We prove the expanding case first. 
An iterate of $u$ under  $\mathcal{E}$ is a point $v$ given
 explicitly as
\[
 v=\begin{cases}
       \displaystyle\frac{1}{2\bar{a}}\left(-\bar{b}u\pm\sqrt{\Delta u^2-4\bar{a} (\bar{\alpha}_1- \bar{\sigma}_1 u)}\right) =: f_\pm(u) &\mbox{if~} \bar{a}\not=0 \\
       \displaystyle\frac{-\bar{\alpha}_1+\bar{\sigma}_1 u-\bar{c}u^2}{\bar{b} u}=: f_0(u) &\mbox{if~} \bar{a}=0, ~\bar{b}\not=0, ~u\not=0.
    \end{cases}
\]
Note that we exclude the case   $\bar{a}=\bar{b}=0$. Because in this case, the relation $\mathcal{E}$ (if exists) is a pair of vertical lines (or a single vertical line if $\bar{\sigma}^2_1=4\bar{\alpha}_1$), thus the slope $m$ is always infinite. If $v=f_0(u)$, then the proof is like that of a $C^1$ interval map on an expanding compact invariant set. Hence, we assume $\bar{a}\not=0$.

Because $\Lambda$ is compact and because we assumed $Df_\pm(u)$ is not infinite when $u\in\Lambda$, there exists $\chi_1>0$ such that $|f_+(u)-f_-(u)|\ge 2\chi_1$ for all $u\in\Lambda$. The continuity of $f_\pm$ guarantees that there exists $\chi_2>0$ such that  $|f_+(u)-f_-(\tilde{u})|\ge \chi_1$ for any $u$, $\tilde{u}\in\Lambda$ with $|u-\tilde{u}|<\chi_2$. Let $\chi=\min\{\chi_1, \chi_2\}$. Now, if $|\xi_0-\tilde\xi_0|\ge \chi_2$, then $|\xi_0-\tilde\xi_0|\ge \chi$. If $|\xi_0-\tilde\xi_0|< \chi_2$, $\xi_n=f_+^n(\xi_0)$ and $\tilde\xi_n=f_+^n(\tilde\xi_0)$ (or $\xi_n=f_-^n(\xi_0)$ and $\tilde\xi_n=f_-^n(\tilde\xi_0)$) for all $n\ge 0$, then the situation is like that of $f_0$ (i.e. $\xi_n$ and $\tilde\xi_n$ are orbit points of a single interval map). Therefore, we  assume  $\xi_N=f_-\circ f_+^{N-1}(\xi_0)$ and $\tilde\xi_N=f_+^N(\tilde\xi_0)$ (or assume $\xi_N=f_+\circ f_-^{N-1}(\xi_0)$ and $\tilde\xi_N=f_-^N(\tilde\xi_0)$) for some $N\ge 1$. If $|\xi_{N-1}-\tilde\xi_{N-1}|\ge \chi_2$ $(\ge \chi)$, then we are done. Otherwise,  $|\xi_{N}-\tilde\xi_{N}|=|f_\mp(\xi_{N-1})-f_\pm(\tilde\xi_{N-1})|\ge \chi_1$ $(\ge \chi)$.

For the contracting case, we consider the inverse iterate of $v$, which is a point $u$ such that $(u,v)\in\mathcal{E}$. Explicitly, 
\[
 u=\begin{cases}
       \displaystyle\frac{1}{2\bar{c}}\left(\bar\sigma_1-\bar{b}v\pm\sqrt{\Delta v^2-4\bar{c} \bar{\alpha}_1+\bar\sigma_1^2- 2\bar{\sigma}_1 \bar{b} v}\right) =: f_\pm(v) &\mbox{if~} \bar{c}\not=0 \\
       \displaystyle\frac{\bar{\alpha}_1+\bar{a}v^2}{\bar\sigma_1-\bar{b}v}=: f_0(v) &\mbox{if~} \bar{c}=0 ~\mbox{and}~ v\not=\bar\sigma_1/\bar{b}.
    \end{cases}
\]
Note that when $\bar{c}=0$ and $v=\bar\sigma_1/\bar{b}$, there does not exist a value $u$ such that point $(u,v)$ is in the relation $\mathcal{E}$ except when $\bar\alpha_1=-\bar{a}\bar\sigma_1^2/\bar{b}^2$.  In the exceptional case, $u=\lim_{v\to\bar\sigma_1/\bar{b}} f_0(v)$. Now, the rest of the proof is the same as that of the expanding case, by considering $f_0$ and $f_\pm$.
\qed

For the relation $\mathcal{E}$, following the terminology of McGehee \cite{McGe1992}, we 
call a subset $\Lambda$ of $\R$ {\it backward complete} if for each $v\in\Lambda$ then there exists $u\in\Lambda$ such that $(u,v)\in\mathcal{E}\cap (\Lambda\times\Lambda)$, and {\it forward complete} 
 if for each $u\in\Lambda$ then there exists $v\in\Lambda$ such that $(u,v)\in\mathcal{E}\cap (\Lambda\times\Lambda)$.
Let $\Lambda$ be a compact subset that is both backward and forward complete. (The existence of such subsets will be addressed in Section \ref{sec:Examples}). 
Define 
\[ 
 \Sigma^\Lambda=\left\{  
                                       (\xi_t)_{t\in\Z}:~(\xi_t, \xi_{t+1})\in \mathcal{E}\cap (\Lambda\times\Lambda)\quad\forall t\in\Z
                              \right\}.
\]

For a map $f$, we call a set $K$   {\it invariant} under $f$ if $f^{-1}(K)=K$.
Denote by $S$ the usual shift operator on $(\R^n)^\Z$ for any integer $n\ge 1$, $S:(z_t)_{t\in\Z}\mapsto (\hat{z}_t)_{t\in\Z}$ with $\hat{z}_t=z_{t+1}$.

\begin{prop} \label{SInvariant}
 $ \Sigma^\Lambda$ is non-empty, compact (with the product topology), and $S$-invariant.
\end{prop}
\proof
 Let $\xi_0\in\Lambda$. Because $\Lambda$ is forward complete, points $\xi_t$ such that $(\xi_t, \xi_{t+1})\in \mathcal{E}\cap (\Lambda\times\Lambda)$ for all $t\ge 0$ exist. Likewise, points $\xi_t$ such that $(\xi_t, \xi_{t+1})\in \mathcal{E}\cap (\Lambda\times\Lambda)$ for all $t\le -1$ exist by the backward completeness. This shows that $\Sigma^\Lambda$ is non-empty. 

It is plain to see that $S$ is a bijection on $\Sigma^\Lambda$. Hence,  $S^{-1}(\Sigma^\Lambda)=\Sigma^\Lambda$.

Let $(\boldsymbol\xi^{(n)})_{n\ge 0}$ be any sequence in $\Sigma^\Lambda$. By the compactness
 of $\Lambda$, there is a convergent subsequence of  $(\xi^{(n)}_0)_{n\ge 0}$. Assume the subsequence, and assume that $\xi^{(n)}_0\to \xi^\infty_0\in \Lambda$ as $n\to\infty$.
Next, consider the sequence $(\xi_1^{(n)})_{n\ge 0}$, and assume  $\xi^{(n)}_1\to \xi^\infty_1\in \Lambda$ as $n\to\infty$. The continuity of the correspondence \eqref{Q} implies $(\xi_0^\infty, \xi_1^\infty)\in\mathcal{E}\cap (\Lambda\times\Lambda)$.  
Similarly, for every integer $t\not=0$ there corresponds  a convergence of points $\xi^{(n)}_t\to \xi^{\infty}_t \in\Lambda$  as $n\to\infty$ such that  
$(\xi^{\infty}_t,  {\xi}^{\infty}_{t+1})\in \mathcal{E}\cap (\Lambda\times\Lambda)$ for all $t\in\Z$. This proves the compactness.  
\qed

\section{Persistence of  AI states} \label{sec:mainthm}

The following is the main theorem of this paper, which says that hyperbolic AI states persist away from the AI limit, becoming genuine orbits of the map $L$. 

\begin{theorem}[Main theorem] \label{mainthm}
 Assume $\Lambda$ is a compact, hyperbolic, backward and forward complete set for the relation $\mathcal{E}(\bar{\alpha}_1, \bar{\sigma}_1, \bar{a}, \bar{c})$. There exists $e_0>0$ such that whenever $\epsilon>0$ and $\|e-e^\dag\|<e_0$ there is a compact hyperbolic  invariant set $\mathcal{A}_e \subset\mathbb{R}^{3}$ for $L$ and a homeomorphism  $\Psi_e:\Sigma^\Lambda\to  \mathcal{A}_e$ such that the following diagram commutes 
 \[
  \begin{matrix}
\Sigma^\Lambda  & \mapright{S} & \Sigma^\Lambda \cr
         \mapldown{\Psi_e}& &\maprdown{\Psi_e} \cr
             \mathcal{A}_e & \mapright{L}& \mathcal{A}_e.  
  \end{matrix}   \tag{$\star$} \label{star}
\]
\end{theorem}
 
The proof of the above theorem will be based on the IFT, applied to $\mathcal{L}(\boldsymbol\xi; e)$ at $e=e^\dag$.

\subsection{Invertibility of $D\mathcal{L}(\boldsymbol\xi; e)$}

The  derivative of $\mathcal{L}$ with respect to $\boldsymbol\xi$ is a linear map 
\[
     D \mathcal{L}(\boldsymbol\xi;e):\Xi\to\Xi, \qquad\boldsymbol\zeta
\mapsto D\mathcal{L}(\boldsymbol\xi;e)\boldsymbol\zeta=\boldsymbol\eta, 
\]
 with
\begin{equation}
  -\epsilon \zeta_{t+1} -\epsilon \sigma \zeta_{t-1}+\epsilon \delta \zeta_{t-2}+2a \xi_t \zeta_t+b\xi_{t-1}\zeta_t+b\xi_{t}\zeta_{t-1}+2c \xi_{t-1}\zeta_{t-1}=\eta_t \label{DFnonhomo}
\end{equation}
for every integer $t$.

\begin{la} \label{mainla}
  If a compact set $\Lambda$ is hyperbolic, backward and forward complete for the relation $\mathcal{E}(\bar\alpha_1, \bar\sigma_1, \bar{a}, \bar{c})$, then $D\mathcal{L}(\boldsymbol\xi; e^\dag)$ is invertible for every $\boldsymbol\xi\in\Sigma^\Lambda$. In addition,
 $\|D\mathcal{L}(\boldsymbol\xi; e^\dag)^{-1}\|_\infty\le 1+C^{-1} (\lambda-1)^{-1}$ with $C$ and $\lambda$ the constants such that \eqref{Clambda+} or \eqref{Clambda-} is satisfied.
\end{la}
\proof
$D\mathcal{L}(\boldsymbol\xi; e^\dag)$ is invertible if and only if
\begin{equation}
(2\bar{a} \xi_t +\bar{b}\xi_{t-1}) \zeta_t-(\bar{\sigma}_1-\bar{b}\xi_{t}-2\bar{c} \xi_{t-1})\zeta_{t-1}=\eta_t   , \quad t\in\Z, \label{non-h} 
\end{equation}
has a unique solution $\boldsymbol\zeta\in\Xi$ for any given  $\boldsymbol\eta\in\Xi$. If $\Lambda$ is  expanding, then $\bar{\sigma}_1-\bar{b}\xi_{t}-2\bar{c}\xi_{t-1}$ cannot be zero for every $t$. The compactness of  $\Lambda$  further implies $\bar{\sigma}_1-\bar{b}\xi_{t}-2\bar{c}\xi_{t-1}$ is uniformly bounded away from zero in $t$. Similarly,    if $\Lambda$ is   contracting, then  $2\bar{a}\xi_{t}+\bar{b}\xi_{t-1}$ is uniformly bounded away from zero in $t$. Thus, we can rescale $\eta_t$ by dividing it by  $\bar{\sigma}_1-\bar{b}\xi_{t}-2\bar{c}\xi_{t-1}$  or $2\bar{a}\xi_{t}+\bar{b}\xi_{t-1}$. And,  instead of \eqref{non-h}, we  consider the following equations 
\begin{equation}
\frac{2\bar{a} \xi_t +\bar{b}\xi_{t-1}}{\bar{\sigma}_1-\bar{b}\xi_{t}-2\bar{c} \xi_{t-1}} \zeta_t-\zeta_{t-1}=\eta_t  \label{etatExpa} 
\end{equation}
for the expanding case, and
\begin{equation}
        \zeta_t-\frac{\bar{\sigma}_1-\bar{b}\xi_{t}-2\bar{c} \xi_{t-1}}{2\bar{a} \xi_t +\bar{b}\xi_{t-1}}\zeta_{t-1}=\eta_t \label{etatCont}
\end{equation}
for the contracting case. In terms of the slope
\[ m_t=m(\xi_t, \xi_{t+1})=\frac{\bar{\sigma}_1-\bar{b}\xi_{t+1}-2\bar{c}\xi_t}{2\bar{a}\xi_{t+1}+\bar{b}\xi_t},
\]
\eqref{etatExpa} means
\begin{equation}
   m_{t-1}^{-1}\zeta_t-\zeta_{t-1}=\eta_t \label{etatm-1}
\end{equation}
and \eqref{etatCont} becomes 
\begin{equation}
  \zeta_t -m_{t-1}\zeta_{t-1}=\eta_t. \label{etatm}
\end{equation}

Equation \eqref{etatm-1} has a solution
\begin{equation}
  \zeta_{t-1}=\zeta_{t-1}(\boldsymbol\eta) =-\eta_t -\sum_{k\ge 0} \left(
                                                                  \prod_{n=0}^k  m_{t-1+n}
                                                     \right)^{-1} \eta_{t+1+k}, \label{zeta+sol} 
\end{equation}
which is bounded for every $t\in\Z$ because 
\begin{eqnarray}
    |\zeta_{t-1}+\eta_t|  \le  \|\boldsymbol\eta\|_{\infty}C^{-1}\sum_{k\ge 0}\lambda^{-k-1} \nonumber \\
\le \|\boldsymbol\eta\|_\infty C^{-1}(\lambda-1)^{-1} \label{norm+}
\end{eqnarray} 
by \eqref{Clambda+}. If there is another solution $(\tilde{\zeta}_t)_{t\in\Z}$, then their difference  $\hat{\zeta}_t=\zeta_t-\tilde{\zeta}_t$ satisfies
\[
m_{t-1}^{-1}\hat{\zeta}_t-\hat{\zeta}_{t-1}=0,
\]
which implies that any non-zero component  $\hat{\zeta}_{t-1}$ can be expressed by 
\[
\hat{\zeta}_{t-1} =    \left( \prod_{n=0}^k  m_{t-1+n}^{-1}
                                                     \right) \hat{\zeta}_{t+k}
\] 
for any $k\ge 0$, and further implies that $\hat{\zeta}_{t+k}$ is  unbounded as $k$ tends to infinity. Therefore, the bounded solution \eqref{zeta+sol} is unique.

Analogously, equation \eqref{etatm} has a solution
\begin{equation}
  \zeta_{t}=\zeta_t(\boldsymbol\eta) =\eta_t +\sum_{k\ge 0} \left(
                                                                  \prod_{n=0}^k  m_{t-1-k+n}
                                                     \right) \eta_{t-1-k}, \label{zeta-sol} 
\end{equation}
which is bounded for every $t\in\Z$ because
\begin{eqnarray}
    |\zeta_{t}-\eta_t|  \le  \|\boldsymbol\eta\|_{\infty}C^{-1}\sum_{k\ge 0}\lambda^{-k-1} \nonumber \\
\le \|\boldsymbol\eta\|_\infty C^{-1}(\lambda-1)^{-1} \label{norm-}
\end{eqnarray} 
by \eqref{Clambda-}. Let $\hat{\zeta}_t=\zeta_t-\tilde{\zeta}_t$ be the difference of $\zeta_t$ with another solution $(\tilde{\zeta}_t)_{t\in\Z}$ (if exists). Then, it   satisfies
\[
 \hat{\zeta}_t-m_{t-1}\hat{\zeta}_{t-1}=0.
\]
Therefore, for some non-zero $\hat{\zeta}_{t}$, we have 
\[
\hat{\zeta}_{t} =    \left( \prod_{n=0}^k  m_{t-1-k+n}      \right) \hat{\zeta}_{t-1-k}
\] 
for any $k\ge 0$. This implies that  $\hat{\zeta}_{t-1-k}$ is  unbounded as $k$ tends to infinity, hence the bounded solution \eqref{zeta-sol} is unique.

By definition, 
\[ \|D\mathcal{L}(\boldsymbol\xi; e^\dag)^{-1}\|_\infty=\sup_{\boldsymbol\eta}\|\boldsymbol\zeta(\boldsymbol\eta)\|_\infty =\sup_{\boldsymbol\eta}\sup_{t\in\Z}|(\zeta_{t}(\boldsymbol\eta)|
\]
with  $\boldsymbol\eta\in\Xi$ and $\|\boldsymbol\eta\|_\infty=1$. By \eqref{norm+} and \eqref{norm-}, and by the triangular inequality, 
\begin{eqnarray*}
 |\zeta_{t}(\boldsymbol{\eta})|
  &\le&        \max\{\eta_t, \eta_{t+1}\}+   \|\boldsymbol\eta\|_\infty C^{-1}(\lambda-1)^{-1} \\
&\le & 1+   C^{-1}(\lambda-1)^{-1}.
\end{eqnarray*}
\qed

\subsection{Proof of Theorem \ref{mainthm}} \label{Proof_main_thm}

Now, we are in a position to prove Theorem \ref{mainthm}. The main idea of proof is similar to that of Theorems 1.3 and A.1 of \cite{Chen2022}.
Let $\boldsymbol\xi^0$ denote a point in $\Sigma^\Lambda$.

\subsubsection*{Continuation from AI states:}

 The linear operator $D\mathcal{L}(\boldsymbol\xi^0; e^\dag)$ is invertible for every $\boldsymbol\xi^0\in \Sigma^\Lambda$  by  Lemma \ref{mainla}.  
Since $\mathcal{L}(\boldsymbol\xi^0; e^\dag)=0$, by virtue of  the IFT   \cite{Zeid1995} there exists $\hat{e}=\hat{e}(\boldsymbol\xi^0)>0$ and a unique $C^1$ function 
  \[
 \theta(\cdot; \boldsymbol\xi^0):\R^5\to \Xi, \qquad e \mapsto   \theta(e;\boldsymbol\xi^0)=(\theta(e; \boldsymbol\xi^0)_t)_{t\in\Z} 
  \]
  such that $\mathcal{L}(\theta(e;\boldsymbol\xi^0); e)=0$ and $\theta(e^\dag; \boldsymbol\xi^0)=\boldsymbol\xi^0$ provided $\|e-e^\dag\|<\hat{e}$.
  Moreover, there exists $\hat{\mu}=\hat{\mu}(\boldsymbol{\xi}^0)>0$ such that if $\mathcal{L}(\boldsymbol{\xi}; e)=0$, $0\le \|e-e^\dag\|<\hat{e}$, 
and $\boldsymbol{\xi}$ belongs to the closed ball $\bar{B}(\boldsymbol{\xi}^0,\hat{\mu})$ of radius $\hat{\mu}$ centered at $\boldsymbol{\xi}^0$ in $\Xi$, then $\boldsymbol{\xi}=\theta(e;\boldsymbol{\xi}^0)$.  
    
   The construction of $\mathcal{L}$ tells that $\theta(e;\boldsymbol\xi^0)$ gives rise to a bounded orbit of the quadratic 3D map $L$ provided $\epsilon >0$.

\subsubsection*{Bijectivity of the continuation:}

For each $\boldsymbol{\xi}^0\in\Sigma^\Lambda$ and $e$, define a map $\mathcal{G}(\cdot; \boldsymbol{\xi}^0, e):\Xi\to \Xi$, $\boldsymbol{\xi}\mapsto \boldsymbol{\xi}-D\mathcal{L}(\boldsymbol{\xi}^0; e^\dag)^{-1}\mathcal{L}(\boldsymbol{\xi},e)$.   By Lemma \ref{mainla}, the norm $\|D\mathcal{L}(\boldsymbol\xi^0; e^\dag)^{-1}\|_\infty$ is bounded above on $\Sigma^\Lambda$. It is clear that 
for any $\gamma>0$ there exist $\mu_0>0$ and $e_1>0$ such that for all $\boldsymbol\xi^0\in\Sigma^\Lambda$ we have $\|\mathcal{L}(\boldsymbol\xi; e)\|_\infty<\gamma$ 
and $\|D\mathcal{L}(\bar{\boldsymbol\eta};  e)-D\mathcal{L}(\boldsymbol\xi^0; e^\dag)\|_\infty<\gamma$ whenever $\|\bar{\boldsymbol\eta}-\boldsymbol\xi^0\|_\infty\le \mu_0$ and $\|e-e^\dag\|\le e_1$. 
 Therefore, we can choose $\mu_0$ and $e_1$ small enough so that
\[ \|D\mathcal{L}(\boldsymbol{\xi}^0; e^\dag)^{-1}\|_\infty ~\|D\mathcal{L}(\bar{\boldsymbol{\eta}}; e)-D\mathcal{L}(\boldsymbol{\xi}^0; e^\dag)\|_\infty \le 1/2\]
and
\begin{equation}
 \|D\mathcal{L}(\boldsymbol{\xi}^0; e^\dag)^{-1}\|_\infty ~\|\mathcal{L}(\boldsymbol{\xi}^0; e)\|_\infty <\mu_0 /2. \label{etadelta}
\end{equation}
Subsequently, for $\boldsymbol{\xi}$, $\boldsymbol{\zeta}\in \bar{B}(\boldsymbol{\xi}^0,\mu_0)$ and $\|e-e^\dag\|\le e_1$, we get 
\begin{eqnarray}
 \lefteqn{ \left\|(\mathcal{G}(\boldsymbol{\xi};\boldsymbol{\xi}^0,e)-\mathcal{G}(\boldsymbol{\zeta};\boldsymbol{\xi}^0,e)\right\|_\infty  } \nonumber \\
&\le &  \sup_{\bar{\boldsymbol{\eta}}\in \{\boldsymbol{\zeta}+t(\boldsymbol{\xi}-\boldsymbol{\zeta}):~ 0\le t\le 1\}}\left\|D\mathcal{L}(\boldsymbol{\xi}^0; e^\dag)^{-1} \left(
                          D\mathcal{L}(\boldsymbol{\xi}^0; e^\dag)-  D\mathcal{L}(\bar{\boldsymbol{\eta}}; e)                                                                                   \right)\right\|_\infty\cdot\|(\boldsymbol{\xi}-\boldsymbol{\zeta})\|_\infty\nonumber\\
&\le& \|\boldsymbol{\xi}-\boldsymbol{\zeta}\|_\infty /2 \label{401}
\end{eqnarray}
and then
\begin{eqnarray}
 \left\|
          \mathcal{G}(\boldsymbol{\xi};
                             \boldsymbol{\xi}^0, e)-\boldsymbol{\xi}^0\right\|_\infty
&\le & \left\| \mathcal{G}(\boldsymbol{\xi}; \boldsymbol{\xi}^0,e)-\mathcal{G}(\boldsymbol{\xi}^0; \boldsymbol{\xi}^0,e)\right\|_\infty+\left\|\mathcal{G}(\boldsymbol{\xi}^0; \boldsymbol{\xi}^0,e)-\boldsymbol{\xi}^0\right\|_\infty\nonumber\\
&\le& \| \boldsymbol{\xi}-\boldsymbol{\xi}^0\|_\infty /2+\|D\mathcal{L}(\boldsymbol{\xi}^0;  e^\dag)^{-1}\|_\infty ~\|\mathcal{L}(\boldsymbol{\xi}^0; e)\|_\infty\label{402}\\
&<&\mu_0. \nonumber
\end{eqnarray}
This implies that $\mathcal{G}(\cdot; \boldsymbol{\xi}^0,e)$ is a contraction map with contraction constant at least $1/2$ on $\bar{B}(\boldsymbol{\xi}^0,\mu_0)$, and $\theta(e;  \boldsymbol{\xi}^0)$ is the unique fixed point in $\bar{B}(\boldsymbol{\xi}^0,\mu_0)$ for $\mathcal{G}(\cdot; \boldsymbol{\xi}^0,e)$
 for any $\boldsymbol{\xi}^0\in\Sigma^\Lambda$
 and $\|e-e^\dag\|\le e_1$. Hence, $\hat{e}(\boldsymbol{\xi}^0)$ and $\hat{\mu}(\boldsymbol{\xi}^0)$ can be chosen so that
$\inf_{\boldsymbol{\xi}^0\in\Sigma^\Lambda}\hat{e}(\boldsymbol{\xi}^0)\ge e_1>0$ and 
 $\inf_{\boldsymbol{\xi}^0\in\Sigma^\Lambda}\hat{\mu}(\boldsymbol{\xi}^0)\ge\mu_0>0$. 
 
We showed in Proposition \ref{Expansive} that  $\Lambda$ is expansive. 
 This implies that $\Sigma^\Lambda$ viewed as a subset of $\Xi$ is uniformly discrete.  
 (A subset $\Sigma$ of $\Xi$ being {\it uniformly discrete} means that there exists $\chi>0$ such that whenever $\boldsymbol{\xi}$ and $\boldsymbol{\eta}$ are distinct points of $\Sigma$, then $\|\boldsymbol{\xi}-\boldsymbol{\eta}\|_\infty >\chi$.)  
 The uniform discreteness further implies that 
  the balls $\bar{B}(\boldsymbol{\xi}^0, \mu_0)$, $\boldsymbol{\xi}^0\in\Sigma^\Lambda$, are disjoint in $\Xi$ provided that $\mu_0$ is sufficiently small. It follows that  the mapping
\[        \Sigma^\Lambda \to  \bigcup_{\boldsymbol\xi^0\in\Sigma^\Lambda}\theta(e; \boldsymbol\xi^0)
   \]
is a bijection provided $\|e-e^\dag\|<e_1$.

\subsubsection*{Conjugacy:}

 Denote by $\pi$ the  following projection and dividing by $\epsilon$ (when $\epsilon \not=0$):  
\begin{equation}
\boldsymbol\xi=(\cdots,\xi_{-1},\xi_0,\xi_1,\cdots)\mapsto (\frac{\xi_0}{\epsilon}, \frac{\xi_{-1}}{\epsilon}, \frac{\xi_{-2}}{\epsilon}) \in\mathbb{R}^3.  \label{projection_pi}
\end{equation}
 Consider the  composition of mappings  
  \[
\boldsymbol\xi^0\ \stackrel{\Phi_e}{\longmapsto}\ \theta(e;\boldsymbol\xi^0)\ \stackrel{\pi}{\longmapsto}~ (\frac{\theta(e;\boldsymbol\xi^0)_0}{\epsilon}, \frac{\theta(e;\boldsymbol\xi^0)_{-1}}{\epsilon}, \frac{\theta(e;\boldsymbol\xi^0)_{-2}}{\epsilon}). 
  \]
  The sequence $(\theta(e;\boldsymbol\xi^0)_t/\epsilon, S^{-1}(\theta(e;\boldsymbol\xi^0))_t/\epsilon, S^{-2}(\theta(e;\boldsymbol\xi^0))_t/\epsilon)_{t\in\Z}$ is  an orbit of $L$ and is uniquely determined by the initial points
 $(\theta(\epsilon;\boldsymbol\xi^0)_0/\epsilon, \theta(\epsilon;\boldsymbol\xi^0)_{-1}/\epsilon,\theta(\epsilon;\boldsymbol\xi^0)_{-2}/\epsilon)$ when $\epsilon>0$. Thus, $\pi$  is bijective.   We showed that  $\Phi_e$ is bijective.  Therefore,  $\pi\circ\Phi_e$ is a bijection. Notice that  $\Sigma^\Lambda$ is  $S$-invariant by
 Proposition \ref{SInvariant}.  Now, 
\begin{eqnarray*}
\lefteqn{\mathcal{L}(\boldsymbol\xi; e)_{t+1}}\\
 &=&  \epsilon^2\alpha -\epsilon \xi_{t+2}-\epsilon \sigma \xi_{t}+\epsilon\delta \xi_{t-1}+a\xi_{t+1}^2+b\xi_{t+1}\xi_{t}+c\xi_{t}^2\\
            &=&\epsilon^2 \alpha -\epsilon S(\boldsymbol\xi)_{t+1}-\epsilon \sigma S(\boldsymbol\xi)_{t-1}+\epsilon \delta S(\boldsymbol\xi)_{t-2}+aS(\boldsymbol\xi)_t^2+bS(\boldsymbol\xi)_{t}S(\boldsymbol\xi)_{t-1}+cS(\boldsymbol\xi)_{t-1}^2 \\
          &=& \mathcal{L}(S (\boldsymbol\xi); e)_t.
\end{eqnarray*}
This means that 
\begin{equation}
  S\circ\mathcal{L}(\boldsymbol{\xi}; e)=\mathcal{L}(S(\boldsymbol{\xi}); e). \label{SLcommute}
  \end{equation} 
  We showed that $\mathcal{L}(\cdot;e)$ has a unique zero at $\theta(e; \boldsymbol{\xi}^0)$ in $\bar{B}(\boldsymbol{\xi}^0;\mu_0)$, so does $\mathcal{L}(S(\cdot);e)$ due to \eqref{SLcommute}. This implies that $\mathcal{L}(\cdot; e)$ has a unique zero at $S\circ\theta(e; \boldsymbol{\xi}^0)$ in $\bar{B}(S(\boldsymbol{\xi}^0),\mu_0)$. Because $\mathcal{L}(\cdot; e)$ has been shown to have a unique zero at $\theta(e; S(\boldsymbol{\xi}^0))$ in $\bar{B}(S(\boldsymbol{\xi}^0), \mu_0)$, it must have $S\circ\theta(e; \boldsymbol{\xi}^0)=\theta(e; S(\boldsymbol{\xi}^0))$ by the uniqueness. Hence, 
  the  following diagram 
\begin{equation}
  \begin{matrix}
        \Sigma^\Lambda& \mapright{S} &   \Sigma^\Lambda \cr
         \mapldown{\pi\circ\Phi_e}& &\maprdown{\pi\circ\Phi_e} \cr
             \mathcal{A}_e& \mapright{L}& \mathcal{A}_e
  \end{matrix} \label{diadiacom}
\end{equation}
commutes for positive $\epsilon$ when $\|e-e^\dag\|<e_1$, where the set $\mathcal{A}_e$ is invariant under $L$ and is defined by 
  \[ \mathcal{A}_e:=\bigcup_{
                                               \boldsymbol\xi^0\in   \Sigma^\Lambda
                                            }                                       \pi\circ\theta(e;\boldsymbol\xi^0).
\]
In other words, $\pi\circ\Phi_e=:\Psi_e$ acts as the  conjugacy.

\subsubsection*{Topological conjugacy}

We have shown in Proposition \ref{SInvariant}
 that  $\Sigma^\Lambda$ is compact within the product topology. The projection $\pi$ defined in \eqref{projection_pi} certainly is continuous for a given $\epsilon$. If $\Phi_e$ is continuous, then  $\pi\circ\Phi_e$ is a continuous bijection from a compact space to a Hausdorff space, thus a homeomorphism. 
Consequently, $\pi\circ\Phi_e$  in  the diagram
 \eqref {diadiacom} serves not only as a conjugacy but a topological conjugacy. The following analysis shows that $\Phi_e$ is continuous.

There exists $0<e_0\le e_1$ such that $\|D\mathcal{L}(\boldsymbol{\xi}^0; e^\dag)^{-1}\|_\infty ~ \|\mathcal{L}(\boldsymbol{\xi}^0; e)\|_\infty<\mu_0/4$ provided $\|e- e^\dag\|\le e_0$.
 Thence, in view of \eqref{401}-\eqref{402}, the ball $\bar{B}(\boldsymbol{\xi}^0,\mu_0/2)$ is mapped into itself under $\mathcal{G}(\cdot; \boldsymbol{\xi}^0,e)$, so $\|\Phi_e (\boldsymbol{\xi}^0)-\boldsymbol{\xi}^0\|_\infty<\mu_0/2$ for all $\boldsymbol{\xi}^0\in\Sigma^\Lambda$ and $0\le \|e-e^\dag\|< e_0$. Given $\boldsymbol{\xi}^0_*= (\xi_{*,i}^0)_{i\in\Z}\in \Sigma^\Lambda$ and  $\left(\boldsymbol{\xi}^{0 (k)}\right)_{k\ge 0}$  a convergent sequence to $\boldsymbol{\xi}^0_*$ in $\Sigma^\Lambda$ (with  the product topology),  suppose that $\left(\Phi_e(\boldsymbol{\xi}^{0 (k)})\right)_{k\ge 0}$ converges, via a subsequence if necessary, 
to a point $\boldsymbol{\xi}^*$ in
the Cartesian product $\overline{U}(\Lambda, \mu_0/2)^\Z$, where $\overline{U}(\Lambda, \mu_0/2)$ is the closed $\mu_0/2$-neighborhood
of $\Lambda$.
 Assume the subsequence. To show $\Phi_e$ is continuous, it suffices to show that 
 $\boldsymbol{\xi}^*=\Phi_e(\boldsymbol{\xi}_*^0)$.  For any $N\in\N$, there is $K\in\N$ such that 
\[ |\xi_i^{0 (k)}-\xi_{*,i}^0|<\mu_0/2\]
for all $k>K$,  $|i|<N$, 
 and then 
\begin{eqnarray*}
|\Phi_e(\boldsymbol{\xi}^{0 (k)})_i-\xi_{*,i}^0| &\le& |\Phi_e(\boldsymbol{\xi}^{0 (k)})_i-\xi_{i}^{0 (k)}|+|\xi_{i}^{0 (k)}-\xi_{*,i}^0|\\
                                                                    &<&\mu_0
\end{eqnarray*}
for $0\le \|e-e^\dag\|<e_0$. Passing to $N\to\infty$, we have 
\begin{equation}
  \|\boldsymbol{\xi}^* -\boldsymbol{\xi}^0_*\|_\infty\le\mu_0.   \label{<=mu0}
\end{equation}
Besides, 
\[ \mathcal{L}(\boldsymbol{\xi}^*; e)=\mathcal{L}(\lim_{k\to\infty}\Phi_e(\boldsymbol{\xi}^{0 (k)}); e)=\lim_{k\to\infty}\mathcal{L}(\Phi_e(\boldsymbol{\xi}^{0 (k)}); e)=0\]
because $\mathcal{L}(\cdot; e)$ is continuous on 
$\overline{U}(\Lambda, \mu_0/2)^\Z$ with respect to the product topology.
By \eqref{<=mu0}, this implies  that $\boldsymbol{\xi}^*$ is a zero of $\mathcal{L}(\cdot; e)$ in $\bar{B}(\boldsymbol{\xi}^0_*,\mu_0)$. Because $\mathcal{L}(\cdot; e)$ on $\bar{B}(\boldsymbol{\xi}^0_*, \mu)$ has a unique zero $\theta(e; \boldsymbol{\xi}_*^0)$, which equals to $\Phi_e(\boldsymbol{\xi}_*^0)$, we conclude that $\boldsymbol{\xi}^*$ must be $\Phi_e(\boldsymbol{\xi}_*^0)$.

\subsubsection*{Hyperbolicity:}

Certainly,  $D\mathcal{L}(\boldsymbol\xi; e)$ is invertible  for any $\boldsymbol\xi \in\Phi_e(\Sigma^\Lambda)$ and $\|e-e^\dag\|<e_0$.
Suppose $\boldsymbol\xi^e =\theta(e;\boldsymbol\xi^0)$, then 
 $\|\boldsymbol\xi^e-\boldsymbol\xi^0\|_\infty$ is uniformly bounded in $\boldsymbol\xi^0$ because $\boldsymbol{\xi}^e\in\bar{B}(\boldsymbol{\xi}^0, \mu_0)$ and $\mu_0$ is independent of $\boldsymbol{\xi}^0$. Therefore, $\|(D\mathcal{L}(\boldsymbol\xi^e; e)-D\mathcal{L}(\boldsymbol\xi^0; e^\dag))\boldsymbol\zeta\|_\infty$ is also uniformly bounded in $\boldsymbol\xi^0$ for a given $\boldsymbol\zeta$, and can be taken to be arbitrarily small by taking a smaller $e$. Since
\[
\|D\mathcal{L}(\boldsymbol\xi^e; e)^{-1}\|_\infty 
= \left(\inf\left\{ \|D\mathcal{L}(\boldsymbol\xi^e; e)\boldsymbol\zeta\|_\infty:~\boldsymbol\zeta\in\Xi \quad\mbox{and}\quad\|\boldsymbol\zeta\|_\infty=1\right\}\right)^{-1},
\]
\[
    \|D\mathcal{L}(\boldsymbol\xi^0; e^\dag)^{-1}\|_\infty= \big(\inf\left\{ \|D\mathcal{L}(\boldsymbol\xi^0; e^\dag)\boldsymbol\zeta\|_\infty:~\boldsymbol\zeta\in\Xi \quad\mbox{and}\quad\|\boldsymbol\zeta\|_\infty=1\right\}\big)^{-1},
\]
and
\[
\|D\mathcal{L}(\boldsymbol\xi^e; e)\boldsymbol\zeta\|_\infty
\ge   \|D\mathcal{L}(\boldsymbol\xi^0; e^\dag)\boldsymbol\zeta\|_\infty  - \|\left(D\mathcal{L}(\boldsymbol\xi^e; e)- D\mathcal{L}(\boldsymbol\xi^0; e^\dag)\right)\boldsymbol\zeta\|_\infty,
\]
we have that $\|D\mathcal{L}(\boldsymbol\xi^e; e)^{-1}\|_\infty$ has a uniform bound in $\boldsymbol\xi^e$ by the uniform boundedness of $\|D\mathcal{L}(\boldsymbol\xi^0; e^\dag )^{-1}\|_\infty$ in $\boldsymbol\xi^0$. 

We recall in Theorem \ref{hypequiv} some equivalent descriptions of hyperbolicity for a diffeomorphism of $\R^k$, $k\ge 2$. 
By Theorem \ref{hypequiv}(iii), the compact invariant set $\mathcal{A}_e$ is hyperbolic for $L$ if and only if for any given $(\eta_t^{(1)},    \eta_t^{(2)},\eta_t^{(3)})_{t\in\Z}\in\l_\infty(\Z,\R^3)$ with  $\|(\eta_t^{(1)},    \eta_t^{(2)},\eta_t^{(3)})_{t\in\Z}\|_\infty=1$ and for any orbit $(x_t, y_t, z_t)_{t\in\Z}$ of $L$ with initial point $(x_0, y_0, z_0)\in\mathcal{A}_e$, the recurrence relation 
\begin{equation}
  \left( 
  \begin{matrix}
        \zeta_{t+1}^{(1)} \\
       \zeta_{t+1}^{(2)} \\
      \zeta_{t+1}^{(3)}
\end{matrix}
\right)-
  \left(
          \begin{matrix}
          2ax_t+by_t & -\sigma+bx_t+2cy_t & \delta \\
1 & 0 &0 \\
 0 &1&0
            \end{matrix} 
  \right)  \left( 
  \begin{matrix}
        \zeta_{t}^{(1)} \\
       \zeta_{t}^{(2)} \\
      \zeta_{t}^{(3)}
\end{matrix}
\right) =  \left( 
  \begin{matrix}
        \eta_t^{(1)} \\
       \eta_t^{(2)} \\
      \eta_t^{(3)}
\end{matrix}
\right), \quad t\in\Z,   \label{Lrecurrence}
\end{equation}
 has a unique solution $(\zeta_t^{(1)},    \zeta_t^{(2)},\zeta_t^{(3)})_{t\in\Z}\in\l_\infty(\Z,\R^3)$ such that  $\|(\zeta_t^{(1)},    \zeta_t^{(2)},\zeta_t^{(3)})_{t\in\Z}\|_{\infty}$ is uniformly bounded in $(x_t, y_t, z_t)_{t\in\Z}$. From the recurrence relation, we have
\[ \zeta_{t+1}^{(2)}= \zeta_{t}^{(1)}+\eta_{t}^{(2)}\quad\mbox{and}\quad  \zeta_{t+1}^{(3)}= \zeta_{t}^{(2)}+\eta_{t}^{(3)}, 
\]
thus  recurrence relation \eqref{Lrecurrence} can be re-written as
\[ 
       \zeta_{t+1}^{(1)} -(2ax_t+by_t)\zeta_t^{(1)}- ( -\sigma+bx_t+2cy_t)  (\zeta_{t-1}^{(1)}+\eta_{t-1}^{(2)}) - \delta ( \zeta_{t-2}^{(1)}+ \eta_{t-2}^{(2)}+\eta_{t-1}^{(3)}) = \eta_{t}^{(1)}, 
\]  
or as
\begin{eqnarray}
    \lefteqn{   \zeta_{t+1}^{(1)} +\sigma\zeta_{t-1}^{(1)}- \delta  \zeta_{t-2}^{(1)}-2ax_t\zeta_t^{(1)}-by_t\zeta_{t}^{(1)} 
-bx_t\zeta_{t-1}^{(1)}
-2cy_t\zeta_{t-1}^{(1)}}    \nonumber\\
  &&\hspace{4cm}  =\eta_{t}^{(1)} +(-\sigma+bx_t+2cy_t)  \eta_{t-1}^{(2)}
+\delta ( \eta_{t-2}^{(2)}+ \eta_{t-1}^{(3)}).  \label{ArriveAt}
\end{eqnarray}
By using $y_t=x_{t-1}$, $z_t=x_{t-2}$, $\xi_t^e=\epsilon x_t$,  letting  $\zeta_t^{(1)}=-\epsilon\zeta_t$, and 
\[
  \epsilon \eta_{t}^{(1)} +(-\epsilon \sigma+b\xi^e_t+2c\xi^e_{t-1})  \eta_{t-1}^{(2)}
+\epsilon \delta ( \eta_{t-2}^{(2)}+ \eta_{t-1}^{(3)})=\epsilon\eta_t,
\]
equation \eqref{ArriveAt} becomes
\[
  -\epsilon \zeta_{t+1} -\epsilon \sigma \zeta_{t-1}+\epsilon \delta \zeta_{t-2}+2a \xi^e_t \zeta_t+b\xi^e_{t-1}\zeta_t+b\xi^e_{t}\zeta_{t-1}+2c \xi^e_{t-1}\zeta_{t-1}=\eta_t, 
\]
which is the same as \eqref{DFnonhomo}, i.e. $D\mathcal{L}(\boldsymbol\xi^e; e)\boldsymbol\zeta=\boldsymbol\eta$.
Hence, by the uniform boundedness of $\|D\mathcal{L}(\boldsymbol\xi^e;e)^{-1}\|_\infty$ in $\boldsymbol\xi^e$ just showed in the last paragraph, we  conclude that \eqref{Lrecurrence} has a unique solution and the solution, which depends on $(x_t, y_t, z_t)_{t\in\Z}$, is uniformly bounded in $(x_t, y_t, z_t)_{t\in\Z}$. 

The proof of Theorem \ref{mainthm} is complete.
\qed

\section{Applications} \label{sec:Examples}
 
In this section, we show that the AI states appearing in \cite{HM2022, HM2024, JLM2008} for given parameters form  compact, hyperbolic, backward and forward complete subsets for  quadratic relations. Hence, the persistence of these AI states to genuine orbits of the 3D map $L$ follows from Theorem \ref{mainthm}.

\subsection{A deterministic one-dimensional map} \label{subsec:JLM}

\begin{figure}[htb]
	\begin{center}
	   \includegraphics[width=0.38\textwidth]{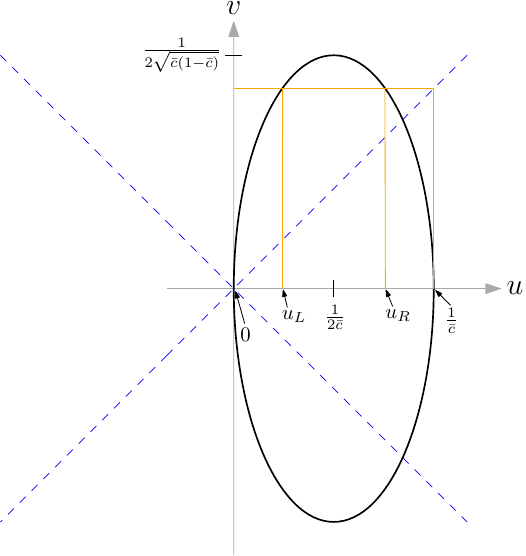}
   \end{center}
	\caption{ A schematic illustration of $\mathcal{E}(0,1, 1-\bar{c}, \bar{c})$ for $0.8<\bar{c}<1$: The relation is depicted in black, and the diagonals $v=\pm u$ are dashed blue. The set $\Lambda$ is a Cantor set and equals $\bigcap_{n=0}^\infty f_+^{-n}([0,1/\bar{c}]\setminus (u_L, u_R))$, where $\{u_L, u_R\}$ is the pre-image of $1/\bar{c}$ under $f_+$.} 
	\label{fig:SomeRelations}
\end{figure}

When $\bar{\alpha}_1=0$, $\bar{\sigma}_1=1$,  $b=O(\epsilon)$ (or simply $\bar{b}=0$) and $0<\bar{c}<1$, we obtain
\begin{equation}
\mathcal{L}(\boldsymbol\xi; e^\dag)_t= (1-\bar{c}) \xi_{t}^2+\bar{c}\xi_{t-1}^2- \xi_{t-1}. \label{eqJLM08}
\end{equation}
 The zeros of  \eqref{eqJLM08} give rise to a relation $\mathcal{E}=\mathcal{E}(0,1,1-\bar{c},\bar{c})$, which is an ellipse (see Figure \ref{fig:SomeRelations}), with $\Delta=-4(1-\bar{c})\bar{c}$.  
 It is plain to see that if $(u,v)\in\mathcal{E}$ then $0\le u\le 1/\bar{c}$. Thus if $v<0$ or $v>1/\bar{c}$ then $(v,w)\not\in\mathcal{E}$ for any $w\in\R$. Consequently, the orbit of the quadratic correspondence \eqref{Q}  must locate on the upper semi-ellipse, and then reduces to that of a one-dimensional map $f_+: \xi_{t-1}\mapsto \xi_t$ given by a branch of the form
\[
    \xi_t=\sqrt{\frac{\xi_{t-1}(1-\bar{c}\xi_{t-1})}{1-\bar{c}}}      
\]
on the domain $[0,1/\bar{c}]$. This is the branch studied in \cite{JLM2008}, where the AI states are orbits of the unimodal map $f_+$.  The map has an escaping interval when $\bar{c}\in (4/5,1)$, and has  positive topological entropy $h_{\rm top}(f_+)$  when $0.791\lessapprox \bar{c}<1$. If $h_{\rm top}(f_+)>0$, then for any $\beta>0$ the map $f_+$ possesses a compact expanding invariant set $\Lambda$   disjoint from the critical point of the  map  
such that $h_{\rm top}(f_+|_\Lambda)>h_{\rm top}(f_+)-\beta$ (see \cite{JLM2008, LM2010}). 

Clearly, the set $\Lambda$ is 
both backward and forward complete because   if $v\in\Lambda$ then there are exactly two $u\in\Lambda$ such that $(u,v)\in\mathcal{E}\cap (\Lambda\times\Lambda)$, and 
 if $u\in\Lambda$ then there is exactly one $v\in\Lambda$ such that $(u,v)\in\mathcal{E}\cap (\Lambda\times\Lambda)$. Moreover, $\Sigma^\Lambda$ is homeomorphic to the inverse limit space\footnote{
 The inverse limit space is defined by $\underleftarrow{\lim} (\Lambda, f_+):=
 \{   ( \ldots, \xi_{-2}, \xi_{-1}, \xi_0
       )|~f_+(\xi_{t-1})=\xi_t\in\Lambda ~\forall t\le 0
  \}$.} $\underleftarrow{\lim} (\Lambda ,f_+)$ of $\Lambda$ for $f_+$.

$\Lambda$ being expanding for $f_+$ precisely means that  there are $C>0$, $\lambda>1$ such that
\[ |Df_+^k(u)|\ge C\lambda^k
\]
 for all $u\in\Lambda$ and $k\ge 1$. And $\Lambda$ turns out to be  expanding for the relation $\mathcal{E}$ because  the above inequality and the one \eqref{Clambda+} are the same  in this case.  Hence, by applying Theorem \ref{mainthm}, there is a compact hyperbolic invariant set $\mathcal{A}_e$ for $L$ such that the diagram \eqref{star} commutes provided
that $\epsilon>0$ and that $e=(\epsilon, \alpha_1, \sigma_1, 0, a, c)$ is sufficiently close to $(0,0, 1, 0, 1-\bar{c}, \bar{c})$ with $0.791\lessapprox\bar{c}<1$ and $a+b+c=1$.

\subsection{A pair of one-dimensional maps} \label{subsec:HM}

When $\bar{\alpha}_1=-1$ and  $\bar{\sigma}_1=r$ for some $r\in\R$, we have 
\begin{equation}
 \mathcal{L}(\boldsymbol\xi; e^\dag)_t=-1 -r \xi_{t-1}+\bar{a}\xi_t^2+\bar{b}\xi_{t}\xi_{t-1}+\bar{c}\xi_{t-1}^2.  \label{HM22}
\end{equation}
 The zeros of  \eqref{HM22} result in  a relation $\mathcal{E}(-1,r,\bar{a}, \bar{c})$. 
When the discriminant $\Delta=-4\bar{a}\bar{c}=\bar{a}r^2$, the relation is a degenerate conic: the ellipse, parabola, and hyperbola turn into a single point, a pair of parallel lines, and two intersecting lines, respectively. 
The relation $\mathcal{E}(-1,0, \bar{a}, \bar{c})$, which is reflection symmetric through the origin,  was  considered in  \cite{HM2022}. When  
 $\Delta$ is equal to zero, it  is  a pair of parallel lines. The relation $\mathcal{E}(-1, r, 1-\bar{c}, \bar{c})$, which is symmetric about the horizontal axis,  was considered in \cite{HM2024}. 

Assuming $\bar{a}\not=0$, an AI-orbit is a sequence $\boldsymbol\xi=(\xi_t)_{t\in\Z}$ satisfying
\begin{equation}
  \xi_t=f_{s_t}(\xi_{t-1})=\frac{1}{2\bar{a}}\left(-\bar{b}\xi_{t-1}+s_t\sqrt{\Delta\xi_{t-1}^2+4\bar{a}(r\xi_{t-1}+1)}\right), \quad s_t\in\{-, +\}, \label{eq14}
\end{equation}
for every integer $t$. 
Providing that the radicand of \eqref{eq14} is strictly positive for all $t$, each AI-orbit is associated with a sequence of binary symbols $(\ldots, s_{-1}, s_0, s_1, \ldots)=:{\bf s}\in \{-,+\}^{\Z}$, also with a sequence of branches $(\ldots, f_{s_{-1}},  f_{s_0}, f_{s_{1}}\ldots)$.  

Let $B\subset \R$ be a compact set, which may depend on $r$, $\bar{a}$ and $\bar{c}$. Given a fixed integer $n\ge 1$, let
\[ \mathcal{R}^+_n=\{(r, \bar{a},\bar{c}):~ \exists B,~ f_{\pm}(B)\subset B,~ 0< |Df_{\bf s}^{\circ n} (x)|\le\lambda^{-1}~\forall x\in B~\mbox{and}~\forall {\bf s}\in\{-,+\}^\Z\}
\]
be a subset in the $(r, \bar{a}, \bar{c})$-space for some $\lambda>1$, 
where $f_{\bf s}^{\circ n}:=f_{s_n}\circ\cdots\circ f_{s_2}\circ f_{s_1}$. (See Remark \ref{condition_imposed}.) Given ${\bf s}\in\{-,+\}^\Z$, define a map $\mathcal{F}(\cdot; {\bf s})=(\mathcal{F}_t(\cdot;{\bf s}))_{t\in\Z}:B^\Z\to B^\Z$
with
\[ \mathcal{F}_t(\boldsymbol\xi; {\bf s})=f_{s_t}(\xi_{t-1})
\]
being the $t^{th}$-component of $\mathcal{F}(\boldsymbol\xi; {\bf s})$.
Note that  
\begin{eqnarray*}
 \mathcal{F}_t^n (\boldsymbol\xi; {\bf s})&=&f_{s_t}  \circ f_{s_{t-1}}\circ\cdots\circ f_{s_{t-n+1}}(\xi_{t-n}) \\
 &=&  f_{S^{t-n}({\bf s})}^{\circ n} (\xi_{t-n})
\end{eqnarray*}
is th  $t^{th}$-component of the $n^{th}$-iterate $\mathcal{F}^n(\boldsymbol\xi; {\bf s})=(\mathcal{F}^n_t(\boldsymbol\xi; {\bf s}))_{t\in\Z}$. 
(We use $S$ to denote both shift operators in the spaces $\{-, +\}^\Z$ and $l_\infty(\Z, \R)$ if it does not cause any ambiguity.)
 If we can show 
\[
\|D\mathcal{F}^n(\cdot; {\bf s})\|_\infty\le \lambda^{-1},
\]
 then
\[
 \|\mathcal{F}^n(\boldsymbol\xi; {\bf s})-\mathcal{F}^n(\boldsymbol\eta;{\bf s})\|_\infty\le \|D\mathcal{F}^n(\cdot; {\bf s})\|_\infty~\|\boldsymbol\xi-\boldsymbol\eta\|_\infty  
\]
for any $\boldsymbol\xi$, $\boldsymbol\eta\in B^\Z$,
thus there exists a unique  fixed point $\phi({\bf s})$ of the map $\mathcal{F}(\cdot; {\bf s})$ by the CMT.  

The  proposition below, which says that  
  there is a one-to-one correspondence between each symbol sequence $(s_t)_{t\in\Z}$ and AI state $(\xi_t)_{t\in\Z}$, is due to Hampton and Meiss.  They proved the proposition for $r=0$ case (see Lemma 1 of \cite{HM2022}) and $\bar{a}=1-\bar{c}$ case (see Lemma 1 of \cite{HM2024}). But, their proofs also work for general $(r, \bar{a}, \bar{c})\in\mathcal{R}_n^+$. Thus, we omit our proof.
\begin{prop} \label{HM'sLemma1}
 Given $(r, \bar{a},\bar{c})\in\mathcal{R}_n^+$ and ${\bf s}=(s_t)_{t\in\Z}\in\{-,+\}^{\Z}$, we have 
\begin{enumerate}
\item 
 $\|D\mathcal{F}^n(\cdot; {\bf s})\|_\infty\le \lambda^{-1}$;
\item 
  there exists exactly one  point, which is $\phi({\bf s})$, in $B^{\Z}$ satisfying \eqref{eq14}; 
\item 
 $\phi({\bf s})\not=\phi({\bf s}^\prime)$ if ${\bf s}\not={\bf s}^\prime\in\{-,+\}^\Z$.
\end{enumerate}
\end{prop}

Let  $\widetilde{\Sigma}$ be the set of AI states obtained from the above proposition. 
The following proposition shows that the one-to-one correspondence given by $\phi$ between symbol sequences and AI states is a continuous map within the product topology.
 
\begin{prop} \label{prop:homeo2}
When $(r, \bar{a}, \bar{c})\in \mathcal{R}^+_n$, with the product topology, $S|_{\widetilde\Sigma}$ is topologically conjugate to $S|_{\{-, +\}^\Z}$.
\end{prop}
\proof
The shift-invariance of $\widetilde\Sigma$ and commutativity of the following diagram 
\[
  \begin{matrix}
\{-, +\}^\Z  & \mapright{S} & \{-, +\}^\Z  \cr
         \mapldown{\phi}& &\maprdown{\phi} \cr
          \widetilde\Sigma  & \mapright{S} & \widetilde\Sigma
  \end{matrix}
\]
are  easy to see: For $\phi({\bf s})=(\phi_t({\bf s}))_{t\in\Z}$ belonging to  $\widetilde{\Sigma}$, we have that   $S^{\pm 1}(\phi({\bf s}))=(\phi_{t\pm 1}({\bf s}))_{t\in\Z}=(\phi_t(S^{\pm 1}({\bf s})))_{t\in\Z}$ also belong to $\widetilde{\Sigma}$. (The last equality follows from Proposition \ref{HM'sLemma1}.) It remains to show that  $\phi$ is a homeomorphism. Since $\{-, +\}^\Z$ is compact and $\widetilde{\Sigma}$ is a Hausdorff space, we only need to show that the bijection $\phi$ is  continuous. 

 Note that $\phi({\bf s})$ is a fixed point of $\mathcal{F}^n(\cdot; {\bf s})$ and that $\phi_t({\bf s})$ is equal to the $t^{th}$-component 
$ \mathcal{F}_t^n (\phi({\bf s}); {\bf s})$ 
of $\mathcal{F}^n(\phi({\bf s}); {\bf s})$. Choose an ${\bf s}^\prime$ such that  $s_{t+i}^\prime=s_{t+i}$ for all $-n+1\le i\le 0$, then
\begin{eqnarray}
|  \phi_t({\bf s}^\prime)-\phi_t({\bf s})|
 &=&  |f_{S^{t-n}({\bf s}^\prime)}^{\circ n} (\phi_{t-n}({\bf s}^\prime)) -f_{S^{t-n}({\bf s})}^{\circ n} (\phi_{t-n}({\bf s})) |  \label{sands1}\\
&=&  |f_{S^{t-n}({\bf s})}^{\circ n} (\phi_{t-n}({\bf s}^\prime)) -f_{S^{t-n}({\bf s})}^{\circ n} (\phi_{t-n}({\bf s})) |  \label{sands2}\\
 &\le &\lambda^{-1} |\phi_{t-n}({\bf s}^\prime)-\phi_{t-n}({\bf s})|    \nonumber
\end{eqnarray}
since $f^{\circ n}_{\bf s}(B)\subset B$ and $|Df_{\bf s}^{\circ n}(x)|\le \lambda^{-1}$ for any $x$ in $B$ and ${\bf s}$ in $\{-,+\}^\Z$.

Choose an ${\bf s}^\prime$ such that  $s_{t-n+i}^\prime=s_{t-n+i}$ for all $-n+1\le i\le 0$, 
and repeat \eqref{sands1} and \eqref{sands2},  then
\[
|  \phi_{t-n}({\bf s}^\prime)-\phi_{t-n}({\bf s})|
 \le \lambda^{-1} |\phi_{t-2n}({\bf s}^\prime)-\phi_{t-2n}({\bf s})|.
\]
Similarly, for each positive integer $k$, we have
\[ 
|  \phi_{t-(k-1)n}({\bf s}^\prime)-\phi_{t-(k-1)n}({\bf s})|
 \le \lambda^{-1} |\phi_{t-kn}({\bf s}^\prime)-\phi_{t-kn}({\bf s})|
\]
for such an ${\bf s}^\prime$ that $s_{t-(k-1)n+i}^\prime=s_{t-(k-1)n+i}$ for all $-n+1\le i\le 0$.  As a result, for an ${\bf s}^\prime$ with $s_{t+i}^\prime=s_{t+i}$ for all $-kn+1\le i\le 0$, we obtain
\begin{eqnarray*}
 |  \phi_{t}({\bf s}^\prime)-\phi_{t}({\bf s})|
 &\le &\lambda^{-k} |\phi_{t-kn}({\bf s}^\prime)-\phi_{t-kn}({\bf s})| \\
&\le & \lambda^{-k} ~ |B|, 
\end{eqnarray*}
where $|B|=\max \{x|~x\in B\}-\min\{x|~x\in B\}$. Consequently, for an ${\bf s}^\prime$ with $s_i^\prime=s_i$ for $-kn+1\le i\le n$, we arrive at 
\[
 |  \phi_{t}({\bf s}^\prime)-\phi_{t}({\bf s})| 
\le \begin{cases}
\lambda^{-k-1}~ |B| & \mbox{for}~  t=n \\
     \lambda^{-k}~ |B| & \mbox{for}~ 0\le t \le n-1 \\
         \lambda^{-k+1} ~ |B| & \mbox{for}~ -n \le t \le -1 \\
    ~~~    \vdots & ~~~\vdots \\
 \lambda^{-1} ~ |B| & \mbox{for}~ -(k-1)n \le t \le -(k-2)n-1. \\
     \end{cases}
\]
Hence, by choosing an ${\bf s}^\prime$ with $s_i^\prime=s_i$ for $-kn+1\le i\le kn$ and passing $k$ to infinity, we conclude that   $\phi({\bf s})$ depends continuously on ${\bf s}$ with the product topology.
\qed

\begin{remark} \label{condition_imposed_2} \rm
When $\bar{c}\not=0$, (by considering reversing time)
an AI state is also a sequence $\boldsymbol\xi=(\xi_t)_{t\in\Z}$ satisfying
\begin{equation}
  \xi_{t-1}=g_{s_t}(\xi_{t})=\frac{1}{2\bar{c}}\left(r-\bar{b}\xi_{t}+s_t\sqrt{\Delta\xi_{t}^2+r^2-2r\bar{b}\xi_{t}+4\bar{c}}\right), \quad s_t\in\{-, +\}, \label{eq38}
\end{equation}
for every integer $t$.  
Let 
\[ \mathcal{R}^-_n=\{(r, \bar{a},\bar{c}):~ \exists B, ~g_{\pm}(B)\subset B,~ 0<|Dg_{\bf s}^{\circ n}(x)|\le \lambda^{-1}~\forall x\in B~\mbox{and}~ \forall{\bf s}\in\{-,+\}^Z\},
\]
where  $g_{\bf s}^{\circ n}:=g_{s_{1-n}}\circ\cdots\circ g_{s_{-1}}\circ g_{s_0}$. Replacing \eqref{eq14} by \eqref{eq38}, 
Propositions \ref{HM'sLemma1} and \ref{prop:homeo2} are valid also for parameter region $\mathcal{R}_n^-$ by constructing a contraction map analogous to $\mathcal{F}(\cdot; {\bf s})$.
\end{remark}

\begin{remark} \label{condition_imposed} \rm
The parameter set $\mathcal{R}_n^+$  defined in \cite{HM2022, HM2024} does not require $0<|Df_{\bf s}^{\circ n}(x)|$, thus is less restricted than the one defined in this paper. We require this condition to rule out the possibility of having an AI orbit point with zero  slope on the relation $\mathcal{E}$. An AI orbit cannot be hyperbolic if there is an orbit point with zero slope. For the same reason, we require  $0<|Dg_{\bf s}^{\circ n}(x)|$ when defining $\mathcal{R}_n^-$. Notice that the parameters in Figures 1 and 4 of \cite{HM2022}, 3(b) and 4(b) of \cite{HM2024} belong to $\mathcal{R}_n^+\cup \mathcal{R}_n^-$ defined here. 
\end{remark}

Finding  AI states is not a trivial task. 
 Remarkably, regions for which $(0, \bar{a},\bar{c})\in \mathcal{R}_1^+\cup\mathcal{R}_1^-$ were found in \cite{HM2022}, and $(r, 1-\bar{c},\bar{c})\in \mathcal{R}_n^+\cup\mathcal{R}_n^-$ for $n$ up to $15$ were obtained in \cite{HM2024}. 

Let $\pi_t$ be the projection to the $t^{th}$-component, $\pi_t(\boldsymbol\xi)=\xi_t$. We obtain  $\pi_t(\widetilde{\Sigma})=:\Lambda$. The set $\Lambda$ is compact, and since $\widetilde{\Sigma}$ is shift-invariant, it is independent of $t$. Obviously, 
 $\widetilde{\Sigma}=\Sigma^\Lambda$.
Direct calculation shows $Df_{s_t}(\xi_{t-1})=m(\xi_{t-1}, \xi_t)$ and $Dg_{s_t}(\xi_t)=m(\xi_{t-1}, \xi_t)^{-1}$. 
The condition  that the radicand in \eqref{eq14} or \eqref{eq38} is strictly positive implies respectively that 
$|  2\bar{a}\xi_t + \bar{b}  \xi_{t-1} |$ or
$|  2\bar{c}\xi_{t-1} -r + \bar{b} \xi_t |$ is also strictly positive, and implies further that 
 the slope $m(  \xi_t, \xi_{t+1}) $ is not infinite nor zero. Therefore,  the conditions $0<|Df_{\bf s}^{\circ n}(x)|\le \lambda^{-1}$ for $\mathcal{R}^+_n$ and $0<|Dg_{\bf s}^{\circ n}(x)|\le \lambda^{-1}$ for $\mathcal{R}^-_n$ mean that  $\Lambda$ is either  contracting or expanding, respectively. 
 At last, $\Lambda$ is again  backward as well as forward complete due to the shift-invariance of $\Sigma^\Lambda$.  
 Hence, by applying Theorem \ref{mainthm}, there is a compact hyperbolic invariant set $\mathcal{A}_e$ for $L$ such that the diagram \eqref{star} commutes provided
that $\epsilon>0$ and that $e=(\epsilon, \alpha_1, \sigma_1, 0, a, c)$ is sufficiently close to $(0,-1, r, 0, \bar{a}, \bar{c})$ with $(r, \bar{a}, \bar{c})\in \mathcal{R}_n^+\cup\mathcal{R}_n^-$ and $a+b+c=1$.
Moreover, with the help of Proposition \ref{prop:homeo2}, an immediate result 
is the following.
\begin{cor} \label{cor:HM_horseshoe}
 Given $(r, \bar{a}, \bar{c})\in\mathcal{R}_n^+ \cup \mathcal{R}_n^-$, the set $\mathcal{A}_e$ in Theorem \ref{mainthm} is a Cantor set.
\end{cor}

\section{Discussions} \label{sec:discussion}

In this section, we discuss some related issues.

\subsection{Infinite Jacobian determinant: $\delta_1\not=0$}   \label{subsec:Mira}

In this subsection, we briefly discuss the situation when $\delta_1=\bar{\delta}_1\not= 0$. In  this situation, instead of \eqref{L0delta2=0}, we have 
\[
 \mathcal{L}(\boldsymbol\xi; e^\dag)_t=\bar{\alpha}_1 -\bar{\sigma}_1 \xi_{t-1}+\bar{\delta}_1\xi_{t-2}+\bar{a}\xi_t^2+\bar{b}\xi_{t}\xi_{t-1}+\bar{c}\xi_{t-1}^2
\]
for $e=e^\dag=(0, \bar\alpha_1, \bar\sigma_1, \bar\delta_1, \bar{a}, \bar{c})$. 
As a consequence, a point $\boldsymbol\xi$ in $\Xi$ for which $\mathcal{L}(\boldsymbol\xi,e^\dag)=0$ is equivalent to a (bi-infinite) bounded orbit of the following two-dimensional endomorphism
\begin{equation}
 (u,v)\mapsto (-\frac{\bar{\alpha}_1-\bar{\sigma}_1 u+\bar{a}v^2+\bar{b}vu+\bar{c}u^2}{\bar{\delta}_1}, u).  \label{MiraEnd}
\end{equation}
This means that at the AI limit $\epsilon\to 0$ the diffeomorphism \eqref{map1} does not reduce to the correspondence \eqref{Q} but the inverse of the above endomorphism.
Note that in the special case that $\bar{b}=\bar{c}=0$ and $\bar{\delta}_1=1$, the endomorphism is called the Mira map \cite{GMO2006}.

The mathematical analysis of perturbing an AI state, which means a (bi-infinite) orbit of  \eqref{MiraEnd}, to a genuine orbit of \eqref{map1} requires a different approach and is 
out of the scope of this paper. In a follow-up paper, we plan to investigate this situation.

\begin{remark} \rm
 We may regard the limit $\bar{\delta}_1\to 0$ as an AI limit for the endomorphism \eqref{MiraEnd}. In the limit, the endomorphism reduces to the correspondence \eqref{Q}.
\end{remark}

\subsection{Invertibility of $D\mathcal{L}(\boldsymbol\xi; e)$ revisited}  \label{subsec:equiva}

Inequality \eqref{Clambda+} implies that there exist $k_1\ge 1$ and $\chi>0$ such that 
\begin{equation}
 \prod_{n= 0}^{k_1-1} |m_{t+n}|\ge 1+\chi ~\quad \forall t\in \Z.   \label{1+chi}
\end{equation}
Analogously, inequality \eqref{Clambda-} implies that 
\begin{equation}
 \prod_{n= 0}^{k_1-1} |m_{t+n}|\le 1-\chi ~\quad \forall t\in \Z.  \label{1-chi}
\end{equation}
Making use of inequalities \eqref{1+chi} and \eqref{1-chi}, we can prove Lemma \ref{mainla} in an alternative way, which is a standard technique of showing invertibility of a linear operator (e.g. Theorem 7.3-1 of \cite{Krey1978}).

\paragraph{\it Alternative proof of Lemma \ref{mainla} (with a different uniform bound).} 
For a given $\boldsymbol\eta\in\Xi$, showing that \eqref{etatm-1} has a unique solution $\boldsymbol\zeta\in\Xi$  amounts to showing that 
\[ \zeta_{t-1}=-\eta_t-\sum_{k=0}^{k_1-2}\left(\prod_{n=0}^km_{t-1+n}^{-1}\right)\eta_{t+1+k}+\left(\prod_{n=0}^{k_1-1}m_{t-1+n}^{-1}\right)\zeta_{t-1+k_1}
\]
has a unique solution for $k_1\ge 2$. (If $k_1=1$, then we just consider \eqref{etatm-1} itself.) Let $\hat\eta_{t}=\eta_t+\sum_{k=0}^{k_1-2}\left(\prod_{n=0}^km_{t-1+n}^{-1}\right)\eta_{t+1+k}$ and $\hat{m}_{t-2+k_1}^{-1}=\prod_{n=0}^{k_1-1}m_{t-1+n}^{-1}$, then it amounts to showing that
\begin{equation}
  -\zeta_{t-1}+\hat{m}_{t-2+k_1}^{-1}\zeta_{t-1+k_1}=\hat\eta_t   \label{eq:Le5of21}
\end{equation}
has a unique solution $\boldsymbol\zeta\in\Xi$ for any given $\hat{\boldsymbol\eta}\in\Xi$ for $k_1\ge 1$. (Note that $\hat\eta_t=\eta_t$ if $k_1=1$.) The above equation can be written as $\hat{\mathcal{B}}\boldsymbol\zeta=\hat{\boldsymbol\eta}$  with $\hat{\mathcal{B}}$ a ``two-diagonal" infinite matrix of the form
\[\hat{\mathcal{B}}_{ts}=\begin{cases}
                                          -1 &\mbox{if}~s=t-1 \\
                                    \hat{m}^{-1}_{t-2+k_1} &\mbox{if}~ s=t-1+k_1 \\
                                           0 &\mbox{otherwise}.
                                     \end{cases}
\]
Clearly, $\hat{\mathcal{B}}=\hat{I}+\hat{M}$, with $\hat{I}$ a ``one-diagonal" infinite matrix and all non-zero entries $-1$, and $\hat{M}$ another ``one-diagonal" infinite matrix with non-zero entries $\hat{m}_{t-2+k_1}^{-1}$.  It is also clear that $\hat{I}$ is invertible, and $\|\hat{I}^{-1}\|_\infty=1$. Now, $\|\hat{\mathcal{B}}-\hat{I}\|_\infty=\|\hat{M}\|_\infty \le 1/(1+\chi)<\|\hat{I}\|_\infty$. Therefore, $\hat{\mathcal{B}}$ is invertible (thus the solution $\boldsymbol\zeta$ exists and is unique), with $\|\hat{\mathcal{B}}^{-1}\|_\infty\le (1+\chi)/\chi$.

For the contracting case, 
   \eqref{etatm} having a unique solution is equivalent to  that 
\[ \zeta_{t}=\eta_t+\sum_{k=0}^{k_1-2}\left(\prod_{n=0}^km_{t-1-n}\right)\eta_{t-1-k}+\left(\prod_{n=0}^{k_1-1}m_{t-1-n}\right)\zeta_{t-k_1}
\]
has a unique solution for $k_1\ge 2$. (If $k_1=1$, we just consider \eqref{etatm} itself.)
 Let $\tilde\eta_{t}=\eta_t+\sum_{k=0}^{k_1-2} (\prod_{n=0}^km_{t-1-n})\eta_{t-1-k}$ and $\tilde{m}_{t-k_1}=\prod_{n=0}^{k_1-1}m_{t-1-n}=\prod_{n=0}^{k_1-1}m_{t-k_1+n}$, then the above equation becomes 
\begin{equation}
 - \tilde{m}_{t-k_1}\zeta_{t-k_1}+\zeta_t=\tilde\eta_t.    \label{eq:Le5of21_2}
\end{equation}
(Again, if $k_1=1$, then  $\tilde\eta_t=\eta_t$.) Let $\tilde{\mathcal{B}}$ be a matrix, with 
\[\tilde{\mathcal{B}}_{ts}=\begin{cases}
                                           -  \tilde{m}_{t-k_1} &\mbox{if}~ s=t-k_1 \\
  1 &\mbox{if}~s=t \\
                                           0 &\mbox{otherwise}.
                                     \end{cases}
\]
Following similar proof as that of the expanding case, we arrive at the invertibility of $\tilde{\mathcal{B}}$, with 
  $\|\tilde{\mathcal{B}}^{-1}\|_\infty\le 1/\chi$.
\qed \\

As a matter of fact, inequalities \eqref{Clambda+} and \eqref{1+chi} imply each other (same to \eqref{Clambda-} and \eqref{1-chi}). The following calculation shows the fact for the expanding (respectively, contracting) case: Let 
\[ m_{\rm min}: = \min_{(u,v)\in\mathcal{E}\cap (\Lambda\times\Lambda)}m(u,v),
\qquad m_{\rm max}: = \max_{(u,v)\in\mathcal{E}\cap (\Lambda\times\Lambda)}m(u,v),
\]
 and let $k=lk_1+i$ for some $l\ge 0$ and $0\le i\le k_1-1$. We assume that  $m_{\rm min}\le 1$ for the expanding case and that $m_{\rm max}\ge 1$ for the contracting case. Otherwise, if $m_{\rm min}>1$ then we can set $k_1=1$ and $1+\chi=m_{\rm min}$ in \eqref{1+chi}. Thus, $k=1$, $C=1$, and $\lambda=m_{\rm min}$ in \eqref{Clambda+}.  Similar, if $m_{\rm max}<1$ then we can set $k_1=1$ and $1-\chi=m_{\rm max}$ in \eqref{1-chi}. Thus,  $k=1$, $C=1$, and $\lambda^{-1}=m_{\rm max}$ in \eqref{Clambda-}.
Thence, 
\begin{eqnarray*} 
\lefteqn{
  \prod_{n= 0}^{k-1} |m_{t+n}| }\\ 
&=& \begin{cases} 
  \displaystyle  \prod_{n= 0}^{k_1-1} |m_{t+n}|\cdot\prod^{k_1-1}_{n=0} |m_{t+k_1+n}|\cdot\ldots\cdot\prod^{k_1-1}_{n=0} |m_{t+(l-1)k_1+n}|\cdot\prod_{n= 0}^{i-1} |m_{t+lk_1+n}| & \mbox{if}~ i >0  \\
  \displaystyle \prod_{n= 0}^{k_1-1} |m_{t+n}|\cdot\prod^{k_1-1}_{n=0} |m_{t+k_1+n}|\cdot\ldots\cdot\prod^{k_1-1}_{n=0} |m_{t+(l-1)k_1+n}|  &  \mbox{if}~ i =0
          \end{cases}                \\
&\ge & (1+\chi)^l \cdot m_{\rm min}^i  \qquad  (\mbox{respectively}, ~ \le (1-\chi)^l \cdot m_{\rm max}^i)\\
&=& (1+\chi)^{(lk_1+i)/k_1} \cdot \frac{m_{\rm min}^i}{(1+\chi)^{i/k_1}} \qquad (\mbox{respectively},~ =(1-\chi)^{(lk_1+i)/k_1} \cdot
 \frac{m_{\rm max}^i}{(1-\chi)^{i/k_1}}) \\
&\ge & (1+\chi)^{k/k_1} \cdot \frac{m_{\rm min}^{k_1-1}}{(1+\chi)^{(k_1-1)/k_1}} \qquad (\mbox{respectively}, ~ \le (1-\chi)^{k/k_1} \cdot \frac{m_{\rm max}^{k_1-1}}{(1-\chi)^{(k_1-1)/k_1}}).
\end{eqnarray*}
By letting $C= {m_{\rm min}^{k_1-1}}/{(1+\chi)^{(k_1-1)/k_1}}$ and $\lambda=(1+\chi)^{1/k_1}$ for the expanding case (respectively, $C^{-1}= {m_{\rm max}^{k_1-1}}/{(1-\chi)^{(k_1-1)/k_1}}$ and $\lambda=(1-\chi)^{-1/k_1}$ for the contracting case), we obtain the desired form for \eqref{Clambda+} (respectively, \eqref{Clambda-}).

\begin{remark} \label{Le5of21}
\rm
 The technique of the above alternative proof of Lemma \ref{mainla} has been employed in literature for the IFT approach to the anti-integrability, for example in the proofs of Lemma 5 of \cite{JLM2008} and Lemma 2.1 of \cite{CLP2016}. The advantage of the original proof in Section  \ref{sec:mainthm} is that it gives the bound of the inverse, namely $\|D\mathcal{L}(\boldsymbol{\xi}; e^\dag)^{-1}\|$, in terms of constants $C$ and $\lambda$ coming from the condition \eqref{Clambda+} or \eqref{Clambda-} on defining the hyperbolicity.  However, if the condition on hyperbolicity is given directly or indirectly as \eqref{1+chi} or \eqref{1-chi}, for example the condition $\mathcal{R}_n^-$ or $\mathcal{R}_n^+$ in Subsection \ref{subsec:HM}, the alternative proof is a convenience since it does not need to construct solutions of \eqref{eq:Le5of21} and \eqref{eq:Le5of21_2}.
\end{remark}

\subsection{Hyperbolicity and nondegeneracy}  \label{equiva2}

We call an AI state $\boldsymbol\xi$ {\it nondegenerate} if the derivative $D\mathcal{L}(\boldsymbol\xi; e^\dag)$ is invertible, otherwise it is called {\it degenerate}. Lemma \ref{mainla} says that an AI orbit lying in a compact hyperbolic set is nondegenerate.  By the IFT, a nondegenerate AI orbit persists, as described in Theorem   \ref{mainthm}.

The notions of hyperbolicity and nondegeneracy  are equivalent, as can be seen by combining  Lemma \ref{mainla} with the lemma below.  

\begin{la} \label{mainla_converse}
Let $\Lambda$ be a compact, backward and forward complete set for the relation $\mathcal{E}$ such that $m(u,v)$ is non-zero and finite for every $(u,v)\in\mathcal{E}\cap (\Lambda\times\Lambda)$. 
If $D\mathcal{L}(\boldsymbol\xi; e^\dag)$ is invertible for every $\boldsymbol\xi\in\Sigma^\Lambda$
 with $\|D\mathcal{L}(\boldsymbol\xi;e^\dag)^{-1}\|_\infty$ bounded uniformly in $\boldsymbol\xi$, then $\Lambda$ is hyperbolic.
\end{la}
\proof
For the expanding case, consider an operator $\mathcal{B}(p; \boldsymbol\xi):\Xi\to\Xi$
with 
\[     (\mathcal{B}(p; \boldsymbol\xi)\boldsymbol\zeta)_t =p\cdot m_{t-1}^{-1}\cdot\zeta_t-\zeta_{t-1}  
\]
for some real number $p\ge 1$. Observe that $(\mathcal{B}(1; \boldsymbol\xi)\boldsymbol\zeta)_t$ is identical to the left hand side of \eqref{etatm-1}. The assumption of $\Lambda$ implies that both $\bar\sigma_1-\bar{b}\xi_t-2\bar{c}\xi_{t-1}$ and $2\bar{a}\xi_t+\bar{b}\xi_{t-1}$ are uniformly bounded away from zero. Therefore, the operator 
$\mathcal{B}( 1; \boldsymbol\xi)$ is invertible.  Moreover, there is $K>0$ such that $\|\mathcal{B}(1;\boldsymbol\xi)^{-1}\|_\infty\le K$ for all $\boldsymbol\xi\in\Sigma^\Lambda$ by our assumption. Now, choose $\tilde\lambda$ with 
$1<\tilde\lambda<1+m_{\rm min}/K$. Then, 
\[ (\mathcal{B}(\tilde\lambda; \boldsymbol\xi)\boldsymbol\zeta)_t-(\mathcal{B}(1; \boldsymbol\xi)\boldsymbol\zeta)_t=(\tilde\lambda-1)m_{t-1}^{-1}\zeta_t.
\]
Thus, 
\[ \|\mathcal{B}(\tilde\lambda; \boldsymbol\xi)-\mathcal{B}(1; \boldsymbol\xi)\|_\infty<K^{-1}\le \|\mathcal{B}(1; \boldsymbol\xi)^{-1}\|_\infty^{-1}.
\]
Therefore, $\mathcal{B}(\tilde\lambda; \boldsymbol\xi)$ is also invertible, with
\begin{equation}
  \|B(\tilde\lambda; \boldsymbol\xi)^{-1}\|_\infty\le (K^{-1}-(\tilde\lambda-1)m_{\rm min}^{-1})^{-1}. \label{Bz-1}
\end{equation}
 Hence, the solution $\boldsymbol\zeta$ of $\mathcal{B}(\tilde\lambda; \boldsymbol\xi)\boldsymbol\zeta=\boldsymbol\eta$ for any $\boldsymbol\eta\in\Xi$ exists in $\Xi$ and is unique, satisfying 
\[ \zeta_{t-1}=-\eta_t-\sum_{k=0}^N\left(\prod_{n=0}^k\tilde\lambda m_{t-1+n}^{-1}\right)\eta_{t+1+k}+\left(\prod_{n=0}^{N+1}\tilde\lambda m_{t-1+n}^{-1}\right)\zeta_{t+1+N}
\]
for any integer $N\ge 0$. For any fixed $t$, choose $\boldsymbol\eta$ such that $|\eta_t|=1$ for all $t$ and  $\left(\prod_{n=0}^{k} m_{t-1+n}^{-1}\right)\eta_{t+1+k} =\left|\prod_{n=0}^{k} m_{t-1+n}^{-1}\right|$ for all $0\le k \le N$. When $N$ tends to infinity, the finiteness of $\zeta_{t-1}$ leads to a sequence $(C_{t, k})_{k\ge 0}$ decreasing to zero as $k$ tends to infinity such that $\left|\prod_{n=0}^{k} m_{t-1+n}^{-1}\right| \tilde\lambda^{k+1}=C_{t,k}$. This implies that $\prod_{n=0}^{k} |m_{t-1+n}| =  \tilde\lambda^{k+1}/C_{t,k}$. Since the value $|\zeta_{t-1}|$ is bounded by the right hand side of \eqref{Bz-1} for all $t$, the value of
 $C_{t,k}$ is bounded above by a positive constant, say $C$, independent of $t$ and $k$. Hence,  $\prod_{n=0}^{k} |m_{t-1+n}| \ge   \tilde\lambda^{k+1}/C$ for all $t\in\Z$ and  $k\ge 0$, as desired.

For the contracting case, let  $\mathcal{B}(p; \boldsymbol\xi)$ be
\[     (\mathcal{B}(p; \boldsymbol\xi)\boldsymbol\zeta)_t =\zeta_t-p\cdot m_{t-1}\cdot\zeta_{t-1}.  
\]
Analogously, there is $K>0$ such that $\|B(1; \boldsymbol\xi)^{-1}\|_\infty\le K$.
Choose $\tilde\lambda$ with 
$1<\tilde\lambda<1+m_{\rm max}^{-1}/K$. Then, 
 $\mathcal{B}(\tilde\lambda; \boldsymbol\xi)$ is  invertible, with
\[
  \|B(\tilde\lambda; \boldsymbol\xi)^{-1}\|_\infty\le (K^{-1}-(\tilde\lambda-1)m_{\rm max})^{-1}.
\]
The unique solution of $\mathcal{B}(\tilde\lambda; \boldsymbol\xi)\boldsymbol\zeta=\boldsymbol\eta$
 satisfies
\[ \zeta_{t}=\eta_t+\sum_{k=0}^N\left(\prod_{n=0}^k\tilde\lambda m_{t-1-n}\right)\eta_{t-1-k}+\left(\prod_{n=0}^{N+1}\tilde\lambda m_{t-1-n}\right)\zeta_{t-2-N}
\]
for any  $N\ge 0$. Similar to the expanding case, there is a sequence $(C_{t, k})_{k\ge 0}$ decreasing to zero and a constant $C$ such that $C_{t,k}\le C$ for all $t\in\Z$ and $k\ge 0$ and that 
$\left|\prod_{n=0}^{k} m_{t-1-n}\right| \tilde\lambda^{k+1}=C_{t,k}$. Hence,  $\prod_{n=0}^{k} |m_{t-1-k+n}| \le  C \tilde\lambda^{-k-1}$.
\qed \\

\begin{remark} \rm
~
\begin{enumerate}
\item
 In the case where $\Sigma^\Lambda$ is obtained by a fixed branch of the correspondence \eqref{Q}, for example the branch  discussed in Subsection \ref{subsec:JLM}, parts of the results of Lemmas \ref{mainla} and \ref{mainla_converse} have been shown as Theorem 2.4 of \cite{Chen2008}. That theorem shows that the hyperbolicity and nondegeneracy are equivalent notions for $C^1$ maps of $\R$. The combination of Lemmas \ref{mainla} and \ref{mainla_converse} generalize that theorem to quadratic correspondences of $\R$. 
\item
For diffeomorphisms in $\R^k$, $k\ge 2$, we also have the equivalence of hyperbolicity and nondegeneracy, to be described in Theorem \ref{hypequiv} below. (Notice that  we have used Theorem \ref{hypequiv}(iii) in the proof of Theorem \ref{mainthm} for the hyperbolicity of the set $\mathcal{A}_e$.)
\end{enumerate}
\end{remark}

In order to formulate Theorem \ref{hypequiv}, let us recall the definition of hyperbolicity for a diffeomorphism. 
  A sequence of invertible linear maps $A_t:\mathbb{R}^k\to\mathbb{R}^k$, $t\in\mathbb{Z}$ and $k\ge 2$, is said to admit a {\it $(\lambda_s,\lambda_u)$-splitting} if there exist direct sum decompositions $\mathbb{R}^k=E_t^s\oplus E_t^u$ such that $A_tE_t^{s/u}=E_{t+1}^{s/u}$ and there are constants $C_s$, $C_u>0$,  $0<\lambda_s <1< \lambda_u$ such that  $\|A_{t+n-1}\cdots A_{t+1}A_t|_{E^s_t}\|\le C_s \lambda_s^n$ and  $\|A_{t-n+1}^{-1}\cdots A_{t-1}^{-1}A_t^{-1}|_{E_{t+1}^u}\|\le C_u \lambda_u^{-n}$ for some norm $\|\cdot\|$ in $\R^k$, all $t\in\mathbb{Z}$ and  $n\ge 1$.    
  Let $f:\mathbb{R}^k\to\mathbb{R}^k$, $k\ge 2$, be a diffeomorphism. Then, a compact $f$-invariant set $\mathcal{A}$ is called  {\it (uniformly) hyperbolic}  for $f$ if for any $w\in \mathcal{A}$ the sequence of differentials $Df(f^t(w)): T_{f^t(w)}\mathbb{R}^k\to T_{f^{t+1}(w)}\mathbb{R}^k$, $t \in\mathbb{Z}$, admits a $(\lambda_s,\lambda_u)$-splitting and the splitting varies continuously with the point $w$. 

Given a  differentiable map  $f$ of $\mathbb{R}^k$, $k\ge 2$, we introduce a map $F$ on the Banach space $l_\infty (\mathbb{Z},\mathbb{R}^k)$:
\begin{eqnarray*}
  F:l_\infty (\mathbb{Z},\mathbb{R}^k) &\to& l_\infty (\mathbb{Z},\mathbb{R}^k),  \\
  (\ldots, w_{-1}, w_0, w_1, \ldots)={\bf w} &\mapsto& F({\bf w})=\left(F({\bf w})_{t}\right)_{t\in\mathbb{Z}} 
\end{eqnarray*}
with $F({\bf w})_t= w_{t+1}-f(w_t)$. One can see that ${\bf w}$ is a bounded orbit of $f$ if and only if it solves $F({\bf w})=0$.
  Certainly, $F$ is a differentiable map, with its derivative at ${\bf w}$ a bounded linear operator:
\begin{eqnarray*}
  DF({\bf w}):l_\infty (\mathbb{Z},\mathbb{R}^k) &\to&  l_\infty (\mathbb{Z},\mathbb{R}^k), \\
   (\zeta_t)_{t\in\Z}=\boldsymbol{\zeta} &\mapsto&  DF({\bf w}; f)\boldsymbol{\zeta}=\left(\sum_{j\in\mathbb{Z}} D_{w_j}F({\bf w}; f)_t~\zeta_j\right)_{t\in\mathbb{Z}}=(\eta_t)_{t\in\Z}=\boldsymbol{\eta}.
\end{eqnarray*}
A bounded orbit ${\bf w}$ is called {\it nondegenerate} if $DF({\bf w})$ is invertible.

\begin{theorem}[see \cite{HK2002, Lanf1985, Math1968}] \label{hypequiv}
 Let $\mathcal{A}$ be a compact invariant set for a diffeomorphism  $f$ of $\mathbb{R}^k$, $k\ge 2$, and ${\bf w}=(w_t)_{t\in\Z}$ an orbit of $f$ with $w_t=f^t(w_0)$ and $w_0\in \mathcal{A}$.  The following statements are equivalent.
\begin{enumerate}
\item $\mathcal{A}$ is  hyperbolic. 
\item For every orbit  ${\bf w}$, the derivative $DF({\bf w})$ is an invertible linear operator of $l_\infty(\mathbb{Z},\mathbb{R}^k)$ with $\|DF({\bf w})^{-1}\|_\infty$ bounded uniformly in $\mathbf{w}$. 
\item For any given $\boldsymbol{\eta}=(\eta_t)_{t\in\Z}\in l_\infty(\mathbb{Z},\mathbb{R}^k)$ with $\|\boldsymbol{\eta}\|_{\infty}=1$ and any orbit ${\bf w}$, the recurrence relation 
  \[
   \zeta_{t+1}-Df(w_t)\zeta_t=\eta_t,\quad t\in\mathbb{Z}, 
  \]
    has a unique solution $\boldsymbol{\zeta}=(\zeta_t)_{t\in\mathbb{Z}}\in l_\infty(\mathbb{Z},\mathbb{R}^k)$ such that $\|\boldsymbol{\zeta}\|_\infty$ is uniformly bounded in ${\bf w}$. 
\end{enumerate}
\end{theorem}

\subsection{Degenerate AI limit} \label{subsec:degenerate}

All AI states discussed in this paper so far are nondegenerate. Here, we demonstrate an example that for $\bar{\sigma}_1=\bar{a}=\bar{c}=0$ every bounded orbit  $(\xi_t)_{t\in\Z}$ of the correspondence \eqref{Q} is degenerate. (We remark that the relation $\mathcal{E}(\bar\alpha_1, 0,0,0)$ is a rectangular hyperbola for nonzero $\bar\alpha_1$.) It is easy to see that $(\xi_t)_{t\in\Z}$ must be a period-$2$ orbit\footnote{When $\bar\alpha_1<0$, the solutions $\xi_t=\pm\sqrt{-\bar\alpha_1}$ of \eqref{xit+1t} for all $t$ give two period-$1$ solutions.} satisfying 
\begin{equation}
   \xi_{t+1}=-\frac{\bar{\alpha}_1}{\xi_t}~\qquad \forall t\in\Z.  \label{xit+1t}
\end{equation}
For example, $(\ldots, \xi_{-1}, \xi_0, \xi_1, \xi_2, \dots)=(\ldots, -\bar{\alpha}_1, 1, -\bar{\alpha}_1, 1, \ldots)$ is a period-$2$ orbit and $(\ldots, \xi_{-1}, \xi_0, \xi_1, \xi_2, \dots)=(\ldots, -\bar{\alpha}_1/2, 2, -\bar{\alpha}_1/2, 2, \ldots)$ is another.  
(We assume that $\bar{\alpha}_1$ is not equal to zero. Then $\xi_t$ cannot be zero for every $t$.)
The product of slopes 
\[ m(\xi_t, \xi_{t+1})\cdot m(\xi_{t+1}, \xi_{t+2})=\frac{\xi_{t+2}}{\xi_t}=1
\]
for every $t$ tells that the orbit cannot be  expanding nor contracting. Now, $e^\dag=(0,\bar\alpha_1, 0,0,0,0)$.  
If $D\mathcal{L}(\boldsymbol\xi; e^\dag)$ is invertible (thus $\boldsymbol\xi$ is nondegenerate), then 
\begin{equation}
 \xi_{t-1}\zeta_t+\xi_t\zeta_{t-1}=\eta_t, \qquad t\in\Z, \label{per-2sol}
\end{equation}
has a unique  solution $\boldsymbol\zeta\in\Xi$ for any given $\boldsymbol\eta\in\Xi$. But, whenever $\boldsymbol\zeta$ is a solution,  $\hat{\boldsymbol\zeta}$ defined by 
\[
 \hat{\zeta}_t= \begin{cases}
                      \zeta_t+1 & \mbox{when $t$ is even} \\
  \zeta_{t}+\displaystyle\frac{\xi_{t}^2}{\bar{\alpha}_1} & \mbox{when $t$ is odd}
\end{cases}
\] 
is also a (bounded) solution: when $t$ is even,
\begin{eqnarray*}
  \xi_{t-1}\hat{\zeta}_t+\xi_t\hat{\zeta}_{t-1} &=&  \xi_{t-1}(\zeta_t+1)+\xi_t (\zeta_{t-1}+\frac{\xi_{t-1}^2}{\bar{\alpha}_1}) \\
&=&  \xi_{t-1}\zeta_t+\xi_t\zeta_{t-1}+\xi_{t-1}-\frac{\bar{\alpha}_1}{\xi_{t-1}}\frac{\xi_{t-1}^2}{\bar{\alpha}_1} \qquad (\mbox{by using} ~ \eqref{xit+1t}) \\
&=&\eta_t, 
\end{eqnarray*}
and when $t$ is odd,
\begin{eqnarray*}
  \xi_{t-1}\hat{\zeta}_t+\xi_t\hat{\zeta}_{t-1} &=&  \xi_{t-1}(\zeta_t+\frac{\xi_{t}^2}{\bar{\alpha}_1})+\xi_t (\zeta_{t-1}+1) \\
&=&  \xi_{t-1}\zeta_t+\xi_t\zeta_{t-1}-\frac{\bar{\alpha}_1}{\xi_{t-2}}\frac{\xi_{t}^2}{\bar{\alpha}_1}+\xi_t \qquad (\mbox{by using} ~ \eqref{xit+1t}) \\
&=&\eta_t \qquad (\mbox{because}~\xi_t=\xi_{t-2}).
\end{eqnarray*}
Hence, any solution of \eqref{per-2sol} is not unique.
 
Owing to the degeneracy, some of the AI states given by \eqref{xit+1t} may not be continued to genuine solutions of \eqref{Recurrence} for small $\epsilon$. As a matter of fact, 
the proposition below tells that except the period-$1$ solutions any bounded solution of \eqref{xit+1t} cannot be continued to a solution of \eqref{Recurrence} for $0<\epsilon<1$, $\alpha=-1/\epsilon^2$, $\sigma=0$, $\delta=1$, and $a=c=0$ (thus $\bar{\alpha}_1=-1$, $\bar{\sigma}_1=\bar{a}=\bar{c}=0$). Note that with these parameters the corresponding 3D diffeomorphism is $(x,y,z)\mapsto (-1/\epsilon^2+z+xy, x,y)$. For example, the period-$2$ solution $(\ldots, \xi_{-1}, \xi_0, \xi_1, \xi_2, \ldots)=(\ldots, 1/2, 2, 1/2, 2, \ldots)$ cannot be continued, but the solution $(\ldots, \xi_{-1}, \xi_0, \xi_1, \xi_2, \ldots)=(\ldots, -1, -1, -1, -1, \ldots)$ can. 

\begin{prop}
When $0<\epsilon<1$ ,  
the following difference equation
 \[
   -1-\epsilon(\xi_{t+1}-\xi_{t-2})+\xi_t\xi_{t-1}=0   
\]
has no period-$2$ solutions. It has exactly two period-$1$ solutions $\xi_t=\pm 1$. 
 \end{prop}
\proof
The periodicity of period-$2$ requires  $\xi_{t+1} = \xi_{t-1}$ and then 
\begin{eqnarray}
 \xi_{t-1} &=& -\frac{1}{\epsilon}+\xi_{t-2}+\frac{\xi_t \xi_{t-1}}{\epsilon}, \label{t-1} \\
 \xi_t &=&-\frac{1}{\epsilon}+\xi_{t-1}-\left(\frac{1}{\epsilon^2}-\frac{ \xi_{t-2} }
                          {\epsilon} - \frac{\xi_t \xi_{t-1}}   {\epsilon^2} \right)~ \xi_t.  \label{t} 
 \end{eqnarray}
 After substituting $\xi_{t-2}$ by $\xi_t$ into \eqref{t-1}, we get 
 \begin{equation}
  \xi_{t-1}=\frac{\epsilon\xi_t-1}{\epsilon-\xi_t}  \label{t-1=}
  \end{equation} 
   provided 
 $\xi_t\not=\epsilon$. Indeed, if $\xi_t=\epsilon$, then $\epsilon$ has to equal $1$, 
but we are only considering the $\epsilon<1$ case.
 Subsequently, $\xi_t$ can be solved from \eqref{t} to get $\xi_t=\pm 1$, implying that $\xi_t=\xi_{t-1}=\xi_{t-2}$ by \eqref{t-1=}.
\qed

\section{Conclusions}

Here, we briefly summarize what are done in this paper by comparing them with the literature. 

For the 3D maps $L$ near the AI limit, if we limit our study to finite Jacobian case, it can be characterized by five  independent parameters $\epsilon$, $\alpha_1$, $\sigma_1$, $a$ and $c$, as explained in the Introduction. The anti-integrability concerns the limiting dynamics when $\epsilon\to 0$ and the  persistence of AI states for small $\epsilon>0$.  

Two cases 
 ${\alpha}_1=0$, ${\sigma}_1=-1$, $c=1$, and
 ${\alpha}_1=-1$, ${\sigma}_1=0$, $a=1$
 were ever considered respectively in \cite{CLP2016} and  \cite{LM2006}.
 For both cases,  the AI states form a two-sided Bernoulli shift on two symbols, and  persist to small $\epsilon>0$ by the classical anti-integrability theory.  

When the AI states are completely determined by a fixed single branch of the quadratic correspondence \eqref{Q}, their persistence to genuine orbits of $L$ was ever studied  in \cite{JLM2008}, and a commutative diagram similar to \eqref{star} was constructed. However, the hyperbolicity of the  genuine orbits continued was not investigated. The main theorem, Theorem \ref{mainthm}, of  this paper contributes to the hyperbolicity of this case (see Subsection \ref{subsec:JLM}, where $\alpha_1=0$, $\sigma_1=1$, $a=1-c$, $0<c<1$ were considered).

The papers considered with the case $\alpha_1=-1$ and $\sigma_1=r$ by Hampton and Meiss \cite{HM2022, HM2024} were the first two to deal with  both branches of the correspondence. Based on the CMT, they proved the persistence of AI states to genuine orbits of $L$, and a bijection between the sets of symbolic sequences and genuine orbits. Under a slightly more restricted condition on $\mathcal{R}^\pm_n$, this paper contributes to the literature by showing that the bijection is indeed a continuous function thus a topological conjugacy. Our main theorem also shows that these genuine orbits form a hyperbolic horseshoe of the map $L$ (see Subsection \ref{subsec:HM}). 

The principal idea of this paper relies on the IFT, which  is a natural tool to study persistence of AI states.
Hyperbolicity of  AI states  is crucial, especially to the IFT approach. (Notice that the condition imposed on $\mathcal{R}_n^\pm$ ensures a class of AI states found by using the CMT in \cite{HM2022, HM2024} are hyperbolic, as discussed in Subsection \ref{subsec:HM}.)  For an AI state $\boldsymbol{\xi}$, its hyperbolicity is a sufficient condition  with which the operator  $D\mathcal{L}(\boldsymbol{\xi}; e^\dag)$ is invertible (see Lemma \ref{mainla}). 
This sufficient condition was also employed to show the invertibility of a linear operator in \cite{JLM2008}, where the hyperbolicity of an AI state comes from the positivity of an interval map's topological entropy. In this paper, we showed that this sufficient condition  is in fact necessary for AI states that are solutions of the quadratic correspondences \eqref{Q} (see Lemma \ref{mainla_converse}, also an example of non-hyperbolic AI states in Subsection \ref{subsec:degenerate}). 

The IFT, however,  does not help us finding new (i.e. dynamically different or distinct from known parameters) AI states. We do not tackle the existence of AI states in this paper. The main theorem assumes the existence of  AI states, then proves their persistence.  Although the approach of \cite{JLM2008} is also via the IFT, the  AI states  there were obtained by knowledge from dynamics of interval maps. The advantage of using the CMT allows Hampton and Meiss \cite{HM2022, HM2024} to find  quite large parameter regions for which nontrivial AI states exist. 
Finding  AI states is crucial to the AI limit too. It would be interesting to find more new nontrivial AI states.

\section*{Acknowledgment}
Useful conversations with James D. Meiss and Amanda E. Hampton  are gratefully acknowledged. 
They read through an early version of this paper and gave valuable comments.
The author thanks the referees, whose comments greatly improved the structure of the paper. He also thanks Tokai University for hospitality where parts of the paper were written. This work 
 was partially supported by NSTC grant  111-2115-M-001-010.

\end{document}